\documentclass[final,a4paper, 11pt]{article}
\RequirePackage{etex} 

\usepackage[letterpaper,top=2cm,bottom=2cm,left=2cm,right=2cm,marginparwidth=1.75cm]{geometry}

\usepackage{enumitem}
\usepackage{tikz}
\usepackage{bm}
\usepackage{matlab-prettifier}
\usepackage{url}

\usepackage{amsfonts}
\usepackage{amsmath}
\usepackage{amssymb}
\usepackage{mathtools}
\usepackage{authblk}
\usepackage{placeins}
\usepackage{hyperref}
\usepackage{cleveref}
\usepackage{autonum}
\usepackage{algorithm}
\usepackage{algorithmic}

\newtheorem{definition}{Definition}
\numberwithin{definition}{section}
\newtheorem{assumption}{Assumption}
\numberwithin{assumption}{section}

\numberwithin{theorem}{section}

\numberwithin{proposition}{section}

\newcommand{\jump}[1]{[\mkern-1.5mu [#1] \mkern-1.5mu]}
\newcommand{\avg}[1]{\{ \mkern-5mu \{#1 \} \mkern-5mu \}}

\newcommand{\lymph}{\texttt{lymph}}

\newcommand{\red}[1]{\textcolor{black}{#1}}

\title{\lymph: discontinuous poLYtopal methods for Multi-PHysics differential problems}

\author[1]{Paola F. Antonietti}
\affil[1]{MOX-Dipartimento di Matematica, Politecnico di Milano, Piazza Leonardo da Vinci 32, Milan, 20133, Italy}
\author[1]{Stefano Bonetti}
\author[1]{Michele Botti}
\author[1]{Mattia Corti}
\author[1]{Ivan Fumagalli}
\author[1]{Ilario Mazzieri}

\begin{document}

\maketitle

\begin{abstract}
We present the library \lymph{} for the finite element numerical discretization of coupled multi-physics problems. \lymph{} is a Matlab library for the discretization of partial differential equations based on high-order discontinuous Galerkin methods on polytopal grids (PolyDG) for spatial discretization coupled with suitable finite-difference time marching schemes. The objective of the paper is to introduce the library by describing it in terms of installation, input/output data, and code structure, highlighting -- when necessary -- key implementation aspects related to the method. A user guide, proceeding step-by-step in the implementation and solution of a Poisson problem, is also provided. In the last part of the paper, we show the results obtained for several differential problems, namely the Poisson problem, the heat equation, the elastodynamics system, and a multiphysics problem coupling poroelasticity and acoustic equations. Through these examples, we show the convergence properties and highlight some of the main features of the proposed method, i.e. geometric flexibility, high-order accuracy, and robustness with respect to heterogeneous physical parameters.
\end{abstract}

\section{Introduction}
\label{sec:introduction}
The numerical solution of coupled multi-physics problems is of crucial importance nowadays, spanning different computational areas and applications. We find coupled problems in several engineering fields, e.g. in the context of life sciences for the modeling of soft tissues such as the heart or the brain, or in computational geosciences for studying seismicity, greenhouse gas sequestration, or geothermal energy production. The numerical simulation of these problems is challenging due to their complex nature: phenomena with different spatial and/or temporal scales, the interaction of several physical laws, and (possibly moving) objects with different materials and properties. Along with their intrinsic complexity, the accurate approximation of such problems through numerical methods often requires different constraints on geometric details, scale resolution, or local refinement of the computational mesh.

Over the past few years, polygonal and polyhedral (polytopal, for short) meshes have become increasingly popular as a solution for these numerical challenges due to their flexibility in representing intricate, possibly moving, geometries and interfaces, and heterogeneous media. Thus, a particular interest has been devoted to the development of numerical methods that can handle general grids, such as Discontinuous Galerkin (see e.g., \cite{Bassi2012,Antonietti2013,Cangiani2014, Cangiani2016, CangianiDongGeorgoulisHouston_2017}), Virtual Element Method (see e.g.,\cite{Beirao2013, BeiraoBrezziMariniRusso_2013b, BeiraoBrezziMariniRusso_2014, BeiraoBrezziMariniRusso_2016, Beirao2018, Beirao2019}), Hybrid High-Order  (see e.g.,\cite{DiPietro2014, DiPietro2017_MathComp, DiPietro2015}), Hybridizable Discontinuous Galerkin  (see e.g.,\cite{Cockburn2008, Cockburn2009, Cockburn2009_MathComp}), and, more recently,  Staggered Discontinuous Galerkin (see e.g.,\cite{Park.Zhao2018, Zhao2019, Zhao2020}). 
On top of geometric flexibility offered by arbitrarily-shaped elements in describing complicated geometries, polytopal meshes exhibit some distinguished advantages over classical meshes, including for example the process of mesh generation and handling, see e.g.,  \cite{Botsch2002,talischi2012polymesher, Levy2015,Vaxman2017,Livesu2019}. \red{On the other hand mesh agglomeration is also a key feature, that does not have an analogous counterpart in classical finite element methods. Agglomeration finds important applications both in the construction of scalable multilevel algebraic solvers and within adaptive algorithms (see e.g. \cite{Antonietti2022_JCP,Antonietti2024} where machine-learning driven refinement and agglomeration strategies for polytopal meshes have been developed). Moreover, agglomerated grids can be effectively employed to accurately represent complicated domain features, such as complex boundaries and/or  embedded structures or microstructures, without the need for excessively fine grids.} 
In this work, we focus on the high-order discontinuous Galerkin finite element method on polytopal grids (PolyDG), see e.g., \cite{Bassi2012,Antonietti2013, Cangiani2014,CangianiDongGeorgoulisHouston_2017}).
The use of the PolyDG method offers numerous benefits when dealing with coupled problems: \textit{(i)} flexible representation of complex geometries, \textit{(ii)} flexibility in refinement and agglomeration strategies, \textit{(iii)} ability to cope with non-conforming interfaces, \textit{(iv)} robustness concerning heterogeneities of physical properties, \textit{(v)} arbitrary polynomial approximation order. Related to the geometrical flexibility, another attractive aspect concerns the treatment of transmission conditions that are usually localized on sub-regions of the computational domain and must be represented without compromising efficiency.
The PolyDG method possesses some distinguishing features, that make it very appealing for multi-physics differential problems. First, it can be easily combined with agglomeration strategies for adaptivity which usually leads to many small faces per element, since the dimension of local approximation space does not depend on the number of faces. Concerning refinement procedures, its hierarchical basis structure can also be exploited. Second, PolyDG methods ensure very good performance in terms of parallelization and scalability -- especially for high polynomial degrees and higher dimensions. 
Another advantage consists in its dimension-independent formulation, which is particularly useful when moving from 2D problems to 3D ones. 
Examples of PolyDG schemes can be found, e.g., in 
\cite{Antonietti2013, Bassi2012} for elliptic problems, in \cite{Houston2002} for advection-diffusion-reaction problems, in \cite{Cangiani2017} for parabolic problems, in \cite{Botti2020_korn, Botti2021} for poroelasticity, 
and in \cite{Antonietti2018} for elastic wave propagation problems. PolyDG discretizations \red{of} multi-physics coupled problems can be found, e.g.,  in  \cite{Corti2023_M3AS, Corti2023, fumagalli2024polytopal} for brain modeling and in  \cite{Antonietti2021,Antonietti2023_SISC,Bonetti2023} for computational geosciences. \\

This paper aims to introduce \texttt{lymph} (discontinuous poLYthopal methods for Multi-PHysics), an open-source MATLAB library, released under the GNU Lesser General Public License, version 3 (LGPLv3), for the PolyDG approximation of multi-physics problems in two-dimensions. We focus on the two-dimensional case to ensure that the syntax remains clear and user-friendly, particularly when employing MATLAB as the software environment for the library. We point out that there is already some software on the market for the numerical approximation of multi-physics problems (e.g. \texttt{Basix} \cite{Scroggs2022}, \texttt{FEniCS} \cite{Alnaes2015}, \texttt{life\textsuperscript{x}} \cite{Africa2022,africa2023lifex,africa2024lifex}, \texttt{MFEM} \cite{Anderson2021}, \texttt{MOOSE} \cite{Perman2020}) and very few for the solution of problems on polytopal meshes (e.g., \texttt{VEM++} \cite{dassi2023vem}, \texttt{HArDCore} \cite{hardcore}, \texttt{MRST} \cite{Lie2019}). However, \texttt{lymph} presents several features that, to the best of the authors' knowledge, make it unique by coupling all the advantages that come from the use of (possibly agglomerated) polytopal meshes and, more specifically, by discretizing the problem via PolyDG schemes. 
The available implementation hinges on Interior Penalty discontinuous Galerkin formulations, but it can be adapted to other DG schemes, such as the local discontinuous Galerkin or the Bassi-Rebay methods, by suitably modifying the assembly and introducing lifting operators.
Moreover, the library is very flexible in terms of coupling different existing physics and implementation of new ones. 
For time dependent problems, time integration is realized by means of finite-difference time integration schemes (Crank-Nicolson for first-order differential systems, Newmark-$\beta$ for second-order ones). \\

The rest of the article is structured as follows: in \cref{sec:polydg} we provide a brief introduction to the PolyDG method, concerning its assumptions and basic elements. In \cref{sec:library} we describe the library in terms of input/output data and code structure, while in \cref{sec:main_lymph} we focus on the core functionalities at the basis of the PolyDG discretization. Then, in \cref{sec:user_guide} we provide a user guide, proceeding step-by-step in the solution of a Poisson problem and we show some results regarding time-dependent problems, in which the main features of the PolyDG method (e.g. geometric flexibility, high-order accuracy, and robustness concerning heterogeneous media) are exploited.\\



\section{Main ingredients of high-order Polytopal discontinuous Galerkin methods}
\label{sec:polydg}

The purpose of this section is to present the mesh assumptions, the discrete spaces, and some preliminary results. 
We introduce a polygonal subdivision $\mathcal{T}_h$ of the computational domain $\Omega\subset\mathbb{R}^2$.
Next, we define the internal edges as 
the intersection of any two neighboring elements of $\mathcal{T}_h$. We define $\mathcal{F}_I$ to be the set of all internal edges. The boundary edges are collected in the set $\mathcal{F}_B$ which yields a subdivision of $\partial\Omega$.  Accordingly, the set of all the edges is given by $\mathcal{F}_h=\mathcal{F}_B\cup\mathcal{F}_I$. In what follows, we introduce the main assumptions on the mesh $\mathcal{T}_h$ (cf. \cite{Cangiani2014,CangianiDongGeorgoulisHouston_2017,Antonietti2021Review}).
\begin{definition}[Polytopic-regular mesh]
\label{def:unif_regular}
A mesh $\mathcal{T}_h$ is polytopic-regular if for any $ \kappa \in \mathcal{T}_h$, there exist a set of non-overlapping simplices contained in $\kappa$, denoted by $\{S_{\kappa}^F\}_{F \subset \partial \kappa}$, such that, for any face $F \subset \partial \kappa$, the following condition holds: \red{$h_{\kappa} \lesssim |S_{\kappa}^F| \ |F|^{-1}$}, with $h_{\kappa}$ denoting the diameter of $\kappa$ (i.e. the maximum distance between two points belonging to the element) and with $| \cdot |$ denoting the Hausdorff measure.
\end{definition}
In the above definition and the following, the symbol $\lesssim$ is used to denote the inequality $x \leq C  y$ for a positive constant $C$ not depending on the mesh size and the polynomial approximation order.
As a basis for the construction of the PolyDG approximation, we define fully discontinuous polynomial spaces on $\mathcal{T}_h$. 
Given an element-wise constant polynomial degree $\ell:\mathcal{T}_h\to\mathbb{N}_{>0}
$ determining the order of the approximation, the discrete spaces are defined such as
\begin{equation}
    \label{eq:discrete_spaces}
    \begin{aligned}
    V_h^{\ell} &= \left\{ v_h \in L^2(\Omega) : v_h |_{\kappa} \in \mathcal{P}^{\ell_{\kappa}}(\kappa) \ \ \forall \kappa \in \mathcal{T}_h \right\}, \quad \red{\mathbf{V}_h^{\ell} = \left[ V_h^{\ell} \right]^2},
    \end{aligned}
\end{equation}
where, for each $\kappa\in\mathcal{T}_h$, the space $\mathcal{P}^{\ell_{\kappa}}(\kappa)$ is spanned by polynomials of maximum degree  $\ell_{\kappa}=\ell_{|\kappa}$. 
We consider a mesh sequence $\{\mathcal{T}_h\}_{h\to0}$ satisfying the following properties:

\begin{assumption}
\label{ass:mesh_Th1}
The mesh sequence $\{\mathcal{T}_h\}_{h\to0}$ and the polynomial degree $\ell$ are such that
\begin{enumerate}[start=1,label={\bfseries A.\arabic*}]
    \item \label{ass:A1} $\{\mathcal{T}_h\}_{h\to0}$ is uniformly polytopic-regular;
    \item \label{ass:A3} For each $\mathcal{T}_h\in \{\mathcal{T}_h\}_{h\to0}$ and for any pair of neighbouring elements $\kappa^+, \kappa^- \in \mathcal{T}_h$, the following $hp$-local bounded variation properties hold: $h_{\kappa^+} \lesssim h_{\kappa^-} \lesssim h_{\kappa^+}$ and $\ell_{\kappa^+} \lesssim \ell_{\kappa^-} \lesssim \ell_{\kappa^+}$.
\end{enumerate}
\end{assumption}
Finally, we also need to introduce the average and jump operators. We start by defining them on each interior edge $F\in\mathcal{F}_I$ shared by the elements $\kappa^{\pm}$ as in \cite{Arnold2002}:
\begin{equation}
    \label{eq:avg_jump_operators}
    \begin{aligned}
    & \jump{a} = a^+ \mathbf{n^+} + a^- \mathbf{n^-}, \ 
    && \jump{\mathbf{a}} = \mathbf{a}^+ \otimes \mathbf{n^+} + \mathbf{a}^- \otimes \mathbf{n^-}, \ 
    &&\jump{\mathbf{a}}_n = \mathbf{a}^+ \cdot \mathbf{n^+} + \mathbf{a}^- \cdot \mathbf{n^-}, \\ 
    & \avg{a} = \frac{a^+ + a^-}{2}, \
    && \avg{\mathbf{a}} = \frac{\mathbf{a}^+ + \mathbf{a}^-}{2}, \ && \avg{\mathbf{A}} = \frac{\mathbf{A}^+ + \mathbf{A}^-}{2},
    \end{aligned}
\end{equation}
where $\mathbf{a} \otimes \mathbf{n} = \mathbf{a}\mathbf{n}^T$, and $a, \ \mathbf{a}, \ \mathbf{A}$ are (regular enough) scalar-, vector-, and tensor-valued functions, respectively. The notation $(\cdot)^{\pm}$ is used to denote the trace on $F$ taken within the interior of $\kappa^\pm$ and $\mathbf{n}^\pm$ is the outer unit normal vector to $\partial \kappa^\pm$. Accordingly, on boundary faces $F\in\mathcal{F}_B$, we set
$$
 \jump{a} = a \mathbf{n},\ \
\avg{a} = a,\ \
\jump{\mathbf{a}} = \mathbf{a} \otimes \mathbf{n},\ \
\avg{\mathbf{a}} = \mathbf{a},\ \
\jump{\mathbf{a}}_n = \mathbf{a} \cdot \mathbf{n},\ \
\avg{\mathbf{A}} = \mathbf{A}.
$$

The ingredients presented above are crucial for deriving the PolyDG semi-discrete formulation of a given partial differential problem. To obtain the fully-discrete formulation of time-dependent problems, the PolyDG discretization in space is coupled with a suitable time-integration scheme (e.g. Crank-Nicolson for first-order problems, Newmark-$\beta$ for second-order problems).


\section{The lymph library}\label{sec:library}

This section introduces the main structure and the components of \lymph~library. We give a high-level overview of the library by describing its main functionalities and folder structure. 
We postpone a more procedural user guide in \cref{sec:user_guide}, going along the numerical approximation and implementation of a specific differential problem.

\subsection*{Overview and installation}
\begin{figure}
    \centering
    \includegraphics[width=0.9\textwidth]{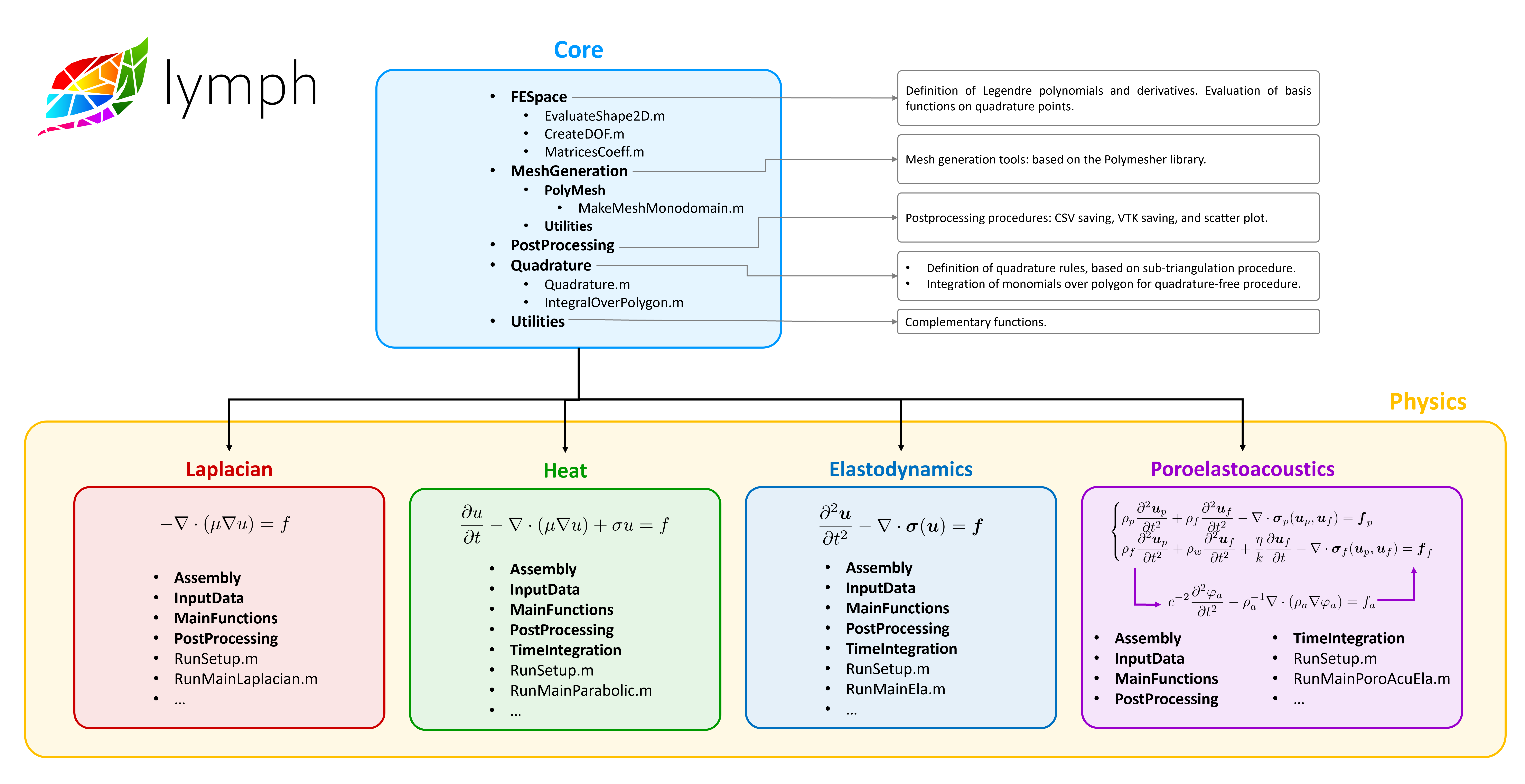}
    \caption{\lymph{} code structure and logo (top left).} 
    \label{fig:scheme}
\end{figure}

\lymph~ is designed to solve either single or multi-physics differential problems employing the PolyDG method, as explained in \cref{sec:polydg}. It has been developed in Matlab version R2022b, but its functionalities have been successfully tested in all versions from R2020b to R2023b including the \texttt{Mapping Toolbox}. The installation of \lymph{} simply consists of downloading the software from the repository \url{https://bitbucket.org/lymph/lymph} into a directory that is accessible to Matlab.
The library is organized into different folders: \texttt{Core} which contains the main routines such as mesh generation, polynomial space construction, and quadrature formulas;
\texttt{Physics}, which contains specific routines for particular problems, \red{organized in subfolders (e.g., \texttt{Laplacian}, \texttt{Heat}, \texttt{Elastodynamics}, \texttt{PoroElastoAcoustics})}. 
A flowchart illustrating the code structure and workflow is given in \cref{fig:scheme}.

\textcolor{black}{
\subsection*{Design and Data Structures: Principles}
The PolyDG discretization is built around several core design principles aimed at flexibility, efficiency, and accuracy in solving PDEs using polygonal meshes. The code leverages specific data structures to efficiently handle mesh generation, finite element computations, and numerical integration to assemble the algebraic system. Below, we detail the key data structures used and the reasons behind their design choices.
\begin{itemize}
\item  \textbf{Mesh Representation}.  Polygonal meshes are represented using a cell array that stores vertex coordinates, element connectivity, and boundary conditions. This structure allows for flexible handling of complex geometries and supports adaptive refinement.
\\
\item \textbf{Finite Element Space (FEspace)}. The FEspace is defined through a structure that includes information about the basis functions, their derivatives, and the degrees of freedom (DOF).
The high-level interface of this centralized structure simplifies the construction of the finite element basis over each mesh element, making the code highly modular and extendable to different element types.
This allows users to easily adapt the FEspace for different polynomial degrees, supporting both low- and high-order approximations.\\
\item \textbf{Quadrature Rules}. Integration over elements and boundaries is handled using a predefined quadrature array that stores the quadrature points and weights. This design allows for precise control over the integration process, which is crucial for accurately assembling the stiffness matrix and load vectors. Using different quadrature rules enables the optimization of computational cost based on the specific problem requirements.\\
\item \textbf{Sparse Matrices for System Assembly}. The assembly of the linear systems, including the stiffness matrix and  the right-hand side, is performed using sparse matrix structures. This is particularly beneficial for large-scale simulations where memory efficiency is crucial. By leveraging sparse storage formats, we minimize memory usage and improve the speed of algebraic solvers applied to the resulting linear systems.\\
\item \textbf{Post-Processing and Visualization}. To facilitate result analysis,  the solution and mesh data  are exported using dedicated structures that can be read by Matlab and Paraview. These structures include solution values at each node and element-level information, which enables seamless post-processing and visualization. This decision was made to ensure compatibility with widely used visualization tools, allowing users to easily interpret and present their simulation results.
\end{itemize}
\subsection*{Customization and User Interaction}
The flexibility of the code is enhanced through user-friendly scripts, such as \texttt{RunSetup.m}, which allow users to configure simulation parameters, select output formats, and define solver settings. This enables a high level of customization without requiring deep modifications to the underlying code. Users can adapt the library to their specific needs, whether they are solving steady-state or time-dependent problems, or validating results through convergence tests.
\subsection*{Documentation and Support}
The full documentation, along with step-by-step tutorials, is provided to guide users through the library's features and applications. This includes practical examples and detailed explanations of the library's setup, aiming to make the onboarding process as smooth as possible for new users. Further information is available at \texttt{https://lymph.bitbucket.io/}, where users can access tutorials and examples to deepen their understanding of the library's capabilities.
}

\section{Common aspects of PolyDG discretization: the \texttt{Core} of \lymph~ }\label{sec:main_lymph}

In this section, we describe the core functionalities of the \lymph~ library, namely, the routines contained in the \texttt{Core} folder that are proper of a PolyDG discretization.  
Moreover, we provide indications on coding strategies that we implemented and can be useful in general finite element software development.

\subsection{Mesh Generation}\label{sec:mesh_generation}
The mesh generation tools included in \lymph{} rely upon \texttt{Polymesher v1.1} \cite{talischi2012polymesher}, an open-source library for the generation of 2D polygonal grids, not directly included in the repository but automatically downloaded at the first run of \lymph{}.
It is important to note that it is possible to use within \lymph{} meshes generated by external software whose output format is compatible with that of \texttt{PolyMesher}. Specifically, agglomerated grids can also be imported and managed.

\subsection{Finite Element Spaces}\label{sec:fespace}
The construction of an algebraic formulation for the discrete problem requires the construction of a basis for the PolyDG space $V_h^\ell$. In the general formulation, the basis is $(\varphi_j)_{j=1}^{N}$, where $N=\dim\{V_h^{\ell}\} = |\mathcal{T}_h|\dim\{P^{\ell}(k)\}$ is the number of degrees of freedom in \texttt{femregion.ndof} with $|\mathcal{T}_h|$ as the number of elements \texttt{femregion.nel} of the partition and $\dim\{P^{\ell}(k)\} = \frac{1}{2}(\ell+1)(\ell+2) =$ \texttt{femregion.nbases}. 
\par
In the \lymph\, library, the basis functions for each physical element $\kappa$ are constructed starting from the Legendre polynomials of order $\ell$ \cite{PoularikasFormulas} in one dimension $\mathcal{L}_i = \mathcal{L}_i(x)$ with $i=0,...,\ell$. Then, by using a tensor product of basis functions in the two directions, we obtain that:
\begin{equation}
    \varphi_i(x,y) = \mathcal{L}_j(x)\mathcal{L}_k(y), \qquad i=1,...,\frac{1}{2}(\ell+1)(\ell+2),\;\mathrm{and}\;j,k=1,...,\ell,\;\mathrm{and}\;j+k\leq\ell. 
\end{equation}
The reference Legendre polynomials $\hat{\mathcal{L}}$ are constructed initially on the square $[-1,1]\times[\-1,1]$; then the final solution is reported on each polygon using an affine transformation to the bounding box of the element, stored in the matrix \texttt{femregion.bbox} \cite{CangianiDongGeorgoulisHouston_2017}.
\par
By construction, the basis of the PolyDG space is modal. Therefore, we cannot associate the coefficients of the linear combination expansion to a physical value at a point in the space. For this reason, the final solution will be reconstructed in the quadrature nodes for visualization purposes. 

\subsection{Evaluation of integrals and use of quadrature formulas}
\label{sec:quadrature}
\red{Any dG discretization requires the computation of volume and boundary integrals. As a general example, we consider the following integral (cf. the Poisson problem in Section~\ref{sec:poisson}) and explain how it is computed in \lymph~. We start by considering the general volume term:}
\begin{equation}\label{eq:A_loc}
    \texttt{A\_loc} \rightarrow (\mu\nabla\varphi_j, \nabla \varphi_i)_{\kappa}, \quad {\rm for} \; i,j =1,..., \frac12(\ell+1)(\ell+2), 
\end{equation}
\red{where $\mu$ is a given element-wise constant function, $\kappa \in \mathcal{T}_h$ is a generic mesh element, and $(\cdot, \cdot)_{\kappa}$ denotes the $L^2(\kappa)$ inner product over $\kappa \in \mathcal{T}_h$. We remark that the syntax $\texttt{A\_loc}$ denotes the volume integral computed for every element $\kappa$ in the mesh $\mathcal{T}_h$; then, for assembling the linear system stemming  from the discretization of our differential equation, we need to combine the local matrices to assemble a global volume matrix $\texttt{A}$.}
The evaluation of these integrals is based on the quadrature-free approach proposed in \cite{AntoniettiHoustonPennesi_18}. The idea of the method is to apply Euler's homogeneous function theorem alongside Stokes' theorem to exactly integrate polynomial functions over polytopal elements. \red{This method allows integrals to be computed without the need for a sub-tessellation of the domain}. \red{In Algorithm~\ref{alg:AssemblyQF}, we report the main steps to compute the values of the volume integrals following the quadrature-free approach (cf. \texttt{lymph/Core/Quadrature} and \texttt{lymph/Laplacian/Assembly} folders for the actual implementation of Algorithm~\ref{alg:AssemblyQF} in \lymph).}
\begin{algorithm}[t]
\color{black}
\caption{\color{black} Assemble volume matrices using the quadrature-free approach of \cite{AntoniettiHoustonPennesi_18}}
\label{alg:AssemblyQF}
\begin{algorithmic}[1]
\STATE Compute the coefficients of the polynomial basis functions
\STATE Compute the reference $2D$-monomial coefficients for integral decomposition
\FOR{every element $\kappa$ in the mesh $\mathcal{T}_h$}
\STATE{Extract the position of the matrix $\texttt{A\_loc}$, cf. \eqref{eq:A_loc}, associated to the current element $\kappa$ w.r.t. the global matrix $\texttt{A}$}
\STATE{Scale the reference $2D$-monomial coefficents according to the shape of the polygon $\kappa$}
\STATE{Compute the bivariate monomial integrals in the polygon $\kappa$}
\STATE{Evaluate the physical parameters}
\STATE{Assemble the local matrix $\texttt{A\_loc}$}
\ENDFOR
\STATE Assemble the global matrix $\texttt{A}$ by combining the local matrices $\texttt{A\_loc}$ along with their positions in the global matrix (cf. Step~4 of this Algorithm).
\end{algorithmic}
\end{algorithm}
\color{black}
\red{We compute the volume integrals relying on the decomposition
$\int_\kappa\varphi_j\varphi_i$ = $\sum_{k,q} a_{kq}^{ij} \int_\kappa x^k y^q$
in terms of bivariate monomials $x^ky^q$.
Once and for all at the beginning, we compute the reference coefficients $\widehat{a}_{kq}^{ij}$ corresponding to the case of $\kappa$ being a rectangle. Then, on each element $\kappa$, we scale the reference coefficients to obtain the current $a_{kq}^{ij}$ and we combine them with the volume integrals of bivariate monomials. For computing the volume integrals we simply need the maximum degree of the monomials in $x$-component, the maximum degree of the monomials in $y$-component, and the geometrical information about the vertices of the polygon (then, of the current element $\kappa \in \mathcal{T}_h$).}
The efficiency of this approach, which does not require the construction of any sub-tessellation of $\kappa$ is discussed in detail in \cite{AntoniettiHoustonPennesi_18}. Notice that the integrals of the monomials are assembled exploiting matrix-vector multiplication, so that the assembly of the local matrices does not need additional \texttt{for} loops that could hinder computational efficiency.
The only $\texttt{for}$ loop needed in the algorithm is the one over mesh elements: this choice allows the code to be flexible w.r.t.~the heterogeneity of physical parameters and ready for using different polynomial degrees $\ell$ on each element ($p$-adaptivity). 
\red{We also remark that the aforementioned algorithm  works also when differential operators (e.g. gradient, divergence) appear in the bilinear forms.}
Another implementation strategy worth pointing out is the construction of each global matrix from the corresponding local ones. We store the local matrices and the corresponding DOF indices in element-indexed arrays, and then pass them to the \texttt{sparse} command to create the global matrices.
This strategy yields a significant gain in computational time w.r.t. preallocating the sparse global matrix and then filling it with local information.

\par
To compute the volume integral of the (possibly non-polynomial) forcing term on the element $\kappa$, we employ a \red{sub-tessellation} of $\kappa$ and use a Gauss-Legendre quadrature rules therein.
More precisely, we consider a sub-tesselation of $\kappa$ made of triangles \texttt{Tria} \red{ and}, on each triangle, we compute the quantities of interest through a quadrature rule with \texttt{femregion.nqn} = \textcolor{black}{$(\ell + 1)^2$} points \cite{zhang_set_2009}.
We denote the corresponding quadrature nodes \texttt{ref\_qNodes\_2D} and weights \texttt{w\_2D}. \red{The quadrature rule is exact for the quantities of interest that involve the trial functions, while the possibility of extending this approach for non-polynomial functions (e.g. forcing terms of generic form) is still under investigation.}
We point out that \red{the quadrature-free formula over each triangle would entail substantially the same computational cost of the chosen quadrature rule, and that} other methods to integrate functions on each (sub) triangle can be applied, also using fewer quadrature nodes \cite{zhang_set_2009,
papanicolopulos_computation_2015,MousaviXiaoSukumar_2010}.
\red{Moreover, we point out that this sub-tessellation strategy can be seen as an interpolatory quadrature formula over the original mesh element $\kappa$, with an interpolation error of the same order of the numerical scheme (or lower): other interpolatory formulas would require the reconstruction of the polynomial interpolation of the forcing term, and since there is no straightforward way to define such interpolation directly on a polygon, this would encompass a computational effort that is comparable to sub-tessellation.}
Although the quadrature-free strategy is preferable for matrix assembly whenever possible, \lymph{} also allows for the implementation of the sub-tessellation strategy. We refer the reader to \texttt{Laplacian/Assembly/MatrixLaplacianST.m} for an example where the use of \texttt{for} loops is minimized, as discussed above.\\

We next discuss the implementation of the face terms. Denoting by $(\cdot, \cdot)_{e}$ the $L^2$ inner product over $e \in \mathcal{F}_I$, we consider the following surface integrals $(\avg{\mu\nabla \varphi_j}, \jump{\varphi_i})_e$ and $(\alpha_e  \jump{\varphi_j}, \jump{\varphi_i})_e$ over an edge $e$ shared by two neighboring elements $\kappa^+$ and $\kappa^-$ in $\mathcal{T}_h$. Using the definitions of $\avg{\cdot}$ and $\jump{\cdot}$ operators, we get:
\begin{itemize}
    \item $\left(\avg{\mu\nabla \varphi_j}, \jump{\varphi_i} \right)_e 
= \frac12((\mu^+\nabla \varphi_j^+  - \mu^-\nabla \varphi_j^-)\cdot \bm n^+, \varphi_i^+)_e + \frac12((\mu^-\nabla \varphi_j^- -  \mu^+\nabla \varphi_j^+) \cdot \bm n^+, \varphi_i^- )_e$,
    \item $\left(\alpha_e  \jump{\varphi_j}, \jump{\varphi_i}\right)_e = (\alpha_e (\varphi_j^+-\varphi_j^-), \varphi_i^+)_e + (\alpha_e (\varphi_j^- - \varphi_j^+), \varphi_i^+)_e $,
\end{itemize}
for $i,j=1,...,\frac12(\ell+1)(\ell+2)$.
When the current element is $\kappa^+$, only the following integrals are 
computed by means of one-dimensional Gauss-Legendre quadrature rules
\begin{equation}\label{eq:IAloc}
       \texttt{IA\_loc}\{\kappa^+\} \rightarrow \frac12(\mu^+\nabla \varphi_j^+ \cdot \bm n^+, \varphi_i^+)_e \quad  {\rm and}  \quad  \texttt{IA\_loc}\{\kappa^+\}  \rightarrow -(\frac12\mu^-\nabla \varphi_j^- \cdot \bm n^+, \varphi_i^+)_e,   
\end{equation}
and
\begin{equation}\label{eq:SAloc}
       \texttt{SA\_loc}\{\kappa^+\}  \rightarrow \alpha_e (\varphi_j^+, \varphi_i^+)_e \quad  {\rm and}  \quad  \texttt{SA\_loc}\{\kappa^+\} \rightarrow -\alpha_e(\varphi_j^-, \varphi_i^+)_e,
\end{equation}
for $i,j =1,..., \frac12(\ell+1)(\ell+2)$. 
\color{black} Boundary integrals for weakly imposing essential boundary conditions are treated in the same way. We omit the discussion for the sake of brevity. In Algorithm~\ref{alg:AssemblyFaces}, we report the main steps to compute the values of the face integrals (cf. \texttt{lymph/Core/Quadrature} and \texttt{lymph/Laplacian/Assembly} folders for the actual implementation of Algorithm~\ref{alg:AssemblyFaces} in \lymph).
\begin{algorithm}[t]
\color{black}
\caption{\color{black} Assemble face matrices}
\label{alg:AssemblyFaces}
\begin{algorithmic}[1]
\STATE Compute reference $1D$ nodes and weights for quadrature on the faces
\FOR{every element $\kappa^+$ in the mesh $\mathcal{T}_h$}
\STATE{Compute the geometric contibution to the penalty coefficient $\alpha_e$ for every edge of the element $\kappa^+$}
\FOR{every edge $e$ of the element $\kappa^+$}
\STATE{Extract the id of the neighboring element $\kappa^-$ (across the current edge $e$)}
\STATE Extract the DOF indexes for assembling the contribution of the current element $\kappa^+$
\STATE{Scale the reference $1D$ nodes and weights according to the edge $e$ (cf. Step~1)}
\STATE{Construct and evaluate the basis functions on the quadrature nodes}
\STATE{Evaluate the physical parameters}
\IF{$e$ is a boundary face}
    \IF{weakly impose boundary conditions on $e$ (e.g. Dirichlet and Robin b.c.)}
        \STATE Assemble the local matrices \texttt{IA\_loc} and \texttt{SA\_loc} (contributions coming from the element $\kappa^+$)
    \ENDIF
\ELSIF{$e$ is an internal face}
    \STATE Assemble the local matrices \texttt{IA\_loc} and \texttt{SA\_loc} (contributions coming from the element $\kappa^+$)
    \STATE Extract the indexes for assembling the contribution of the neighboring element $\kappa^-$
    \STATE Update the local matrices \texttt{IA\_loc} and \texttt{SA\_loc} with the contributions coming from the element $\kappa^-$, cf. \eqref{eq:IAloc} and \eqref{eq:SAloc}, respectively
\ENDIF
\ENDFOR
\ENDFOR
\STATE Assemble the global matrices $\texttt{IA}$, $\texttt{SA}$ by combining the local matrices $\texttt{IA\_loc}$, $\texttt{SA\_loc}$ along with their positions in the global matrix (cf. Step~6, Step~16 in the Algorithm).
\end{algorithmic}
\end{algorithm}
\color{black}

\section{Examples}\label{sec:user_guide}
In the following, we show how to solve differential problems with \lymph~.
In the subsequent discussion, we present four examples of increasing complexity: the Poisson equation, the heat equation, the elastodynamics equation, and a multiphysics system modeling wave propagation in poroelasto-acoustic media. These examples are paradigmatic in illustrating the capabilities of \lymph~ in handling single-physics problems—encompassing both steady-state and time-dependent, as well as scalar and vector-valued differential equations—alongside multiphysics models.

\subsection{The Poisson problem}
\label{sec:poisson}
We start, by considering the following problem in a polygonal domain $\Omega \subset \mathbb{R}^2$:
\begin{equation}\label{eq:poisson}
    \begin{cases}
        - \nabla \cdot (\mu \nabla u) (\bm x) = f(\bm x), & \bm x \in \Omega,\\
        u (\bm x) =  g(\bm x), & \bm x \in \partial \Omega,
    \end{cases}
\end{equation}
where $\mu, f$ and $g$ are given regular functions. We reformulate it using the PolyDG discretization described in \cref{sec:polydg} obtaining:
find $u_h \in V_h^\ell$ s.t. 
\begin{equation}\label{eq:dgPoisson}
    a_{dG}(u_h, v_h) = F(v_h) \quad \forall \; v_h \in  V_h^\ell, 
\end{equation}
where
\begin{equation}
    a_{dG}(u,v) = \sum_{\kappa \in \mathcal{T}_h} (\mu \nabla u, \nabla v)_{\kappa} - \sum_{e \in \mathcal{F}_h}\Big( \big( \avg{\mu\nabla u}, \jump{v} \big)_e + \big( \jump{u}, \avg{\mu \nabla v} \big)_e -  \big( \alpha_e  \jump{u}, \jump{v}\big)_e \Big) 
    \quad \forall u,v \in V_h^\ell,
\end{equation}
with the penalization parameter $\alpha: \mathcal{F}\rightarrow\mathbb{R}_+$ defined as \cite{CangianiDongGeorgoulisHouston_2017}:
\begin{equation}\label{def:penalty}
    \alpha_e(\bm x) = \begin{cases}
        C_\alpha \max_{\kappa \in \{ \kappa^+,\kappa^-\}} \left(\mu_\kappa \frac{\ell_\kappa^2}{h_\kappa} \right), & \bm x \in e, e \in \mathcal{F}_I, e \subset \partial \kappa^+ \cap \partial \kappa^-,  \\
        C_\alpha \mu_\kappa \frac{\ell_\kappa^2}{h_\kappa}, & \bm x \in e, e \in \mathcal{F}_B, e \subset \partial \kappa^+ \cap \partial \Omega,
    \end{cases}
\end{equation}
with $C_\alpha>0$ denoting the penalty coefficient to be properly set, and 
\begin{equation}
    F(v) = \sum_{\kappa \in \mathcal{T}_h} (f,v)_{\kappa} - \sum_{e \in \mathcal{F}_B} \Big( (g, \mu \nabla v)_e -  (\alpha_e g, v)_e \Big) \quad \forall v \in V_h^\ell.
\end{equation}
By introducing a set of basis functions $\{ \varphi_j\}_{j=1}^{N_h}$ for the space $V_h^\ell$
we can write \eqref{eq:dgPoisson} as the following algebraic problem: find $\bm U_h \in \mathbb{R}^{N_h}$ s.t. 
\begin{equation}\label{eq:linear_system_poisson}
    A_{dG} \bm U_h = \bm F,
\end{equation}
with $A_{dG} \in \mathbb{R}^{N_h} \times \mathbb{R}^{N_h}$ defined  for any $i,j = 1,..., N_h$ as $(A_{dG})_{ij} = a_{dG}(\varphi_j,\varphi_i)$
and $\bm F \in \mathbb{R}^{N_h} $ is given by $\bm F_i = F(\varphi_i)$ for all $i = 1,..., N_h$.
The entries of the matrix $A_{dG}$ as well as the right-hand side $\bm F$ in \eqref{eq:linear_system_poisson} are computed as explained in \cref{sec:main_lymph}. For the sake of completeness, we introduce the dG-norm as  $\| u \|_{dG}^2 = \| \sqrt{\mu} \nabla u \|^2_{L^2(\Omega)} + \|\sqrt{\alpha_e} \jump{u} \|^2_{L^2(\mathcal{F}_h)} $ for any $u\in V_h^\ell$, and we recall a well-known convergence result of the PolyDG discretization, \cite[Theorem 36]{CangianiDongGeorgoulisHouston_2017, Antonietti2016Cangiani} from which we have the 
following convergence rates:
\begin{equation}
 \label{eq:error_estimates}
||| u-u_h |||_{dG}^2 \lesssim \sum_{\kappa \in \mathcal{T}_h} \frac{h_\kappa^{2(s_\kappa-1)}}{\ell_\kappa^{2(m_\kappa-1)}} \| u \|_{H^{m_k}(\kappa)}^2,
\quad 
\| u-u_h \|_{L^2(\Omega)}^2 \lesssim \sum_{\kappa \in \mathcal{T}_h} \frac{h_\kappa^{2s_\kappa}}{\ell_\kappa^{2 m_\kappa}} \| u \|_{H^{m_\kappa}(\kappa)}^2,
\end{equation}
with $s_\kappa = \min(\ell_\kappa + 1, m_\kappa)$ for all  $\kappa \in \mathcal{T}_h$, $m_\kappa>0$ denoting the local Sobolev regularity of the solution $u$, and $h_{\kappa}$ the element diameter.

\subsubsection{Verification test}\label{sec:Poisson_verification}
\begin{figure}
    \centering
    \includegraphics[width=0.75\textwidth]{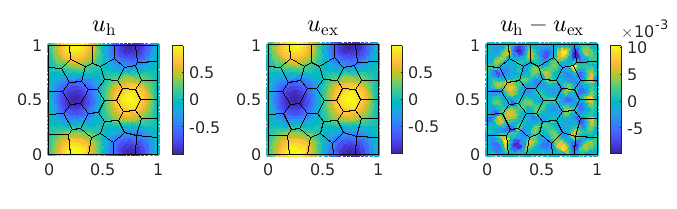}
    \caption{Test case of Section~\ref{sec:poisson}. Left: computed PolyDG solution $u_h$ using a polygonal mesh with $N_{el} = 30$ elements and a polynomial degree $\ell = 3$. Center: analytical solution $u_{ex}$. Right: the difference between numerical and analytical solutions.}
    \label{fig:Poisson_Nel30_l3}
\end{figure}

\begin{figure}
    \centering
    \includegraphics[width=0.45\textwidth]{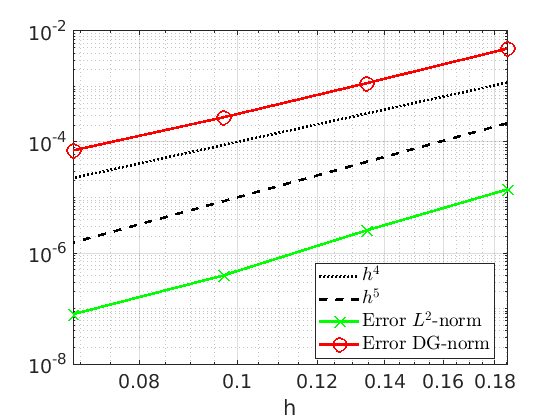}
    \includegraphics[width=0.45\textwidth]{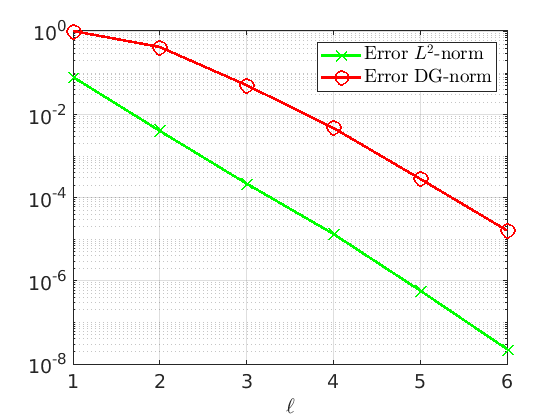}
    \caption{Test case of Section~\ref{sec:poisson}. Left: computed errors $\|u-u_h\|_{dG}$ and $\|u-u_h\|_{L^2(\Omega)}$ as a function of the mesh size $h$ by fixing the polynomial degree $\ell=4$. Right: computed errors $\|u-u_h\|_{dG}$ and $\|u-u_h\|_{L^2(\Omega)}$ as a function of the polynomial degree $\ell=4$ by fixing the number of mesh element $N_{el}=100$.} 
    \label{fig:convergence_laplacian}
\end{figure}

We consider problem \eqref{eq:poisson} in $\Omega = (0,1)^2$ with the following data: $\mu(\boldsymbol{x}) = 1$, \\ $f(\boldsymbol{x}) = 8\pi^2\sin(2\pi x)\cos(2\pi y)$, and $g(\boldsymbol{x}) = \sin(2\pi x)\cos(2\pi y)$, whose exact solution is \\ \noindent $u(\boldsymbol{x}) = \sin(2\pi x)\cos(2\pi y)$. To solve this problem we use the functions contained in \texttt{Laplacian}. We set up these data in \texttt{InputData/DataTestLap.m},  and fix the number of the element mesh $N_{el}=30$, the polynomial approximation degree $\ell_\kappa=3$ for any $\kappa \in \mathcal{T}_h$, and the penalty constant $C_\alpha = 10$ in \eqref{def:penalty}. Next, we run the simulation using the script \texttt{RunMainLaplacian.m},  which calls the main algorithm  \texttt{MainFunctions/MainLaplacian.m}.
 As the output of the run we obtain the plots in Figure~\ref{fig:Poisson_Nel30_l3} showing the
computed solution $u_h$ (left), the analytical solution $u_{ex}$ (center), and the arithmetic difference between the two (right).
Moreover, the output structure \texttt{Error} contains the following fields: \texttt{Nel = 30} (number of mesh elements), \texttt{h = 0.3235} (mesh size), \texttt{p = 3} (polynomial approximation degree), \texttt{L2 = 0.0027} ($L^2$-norm of the error), \texttt{dG = 0.3349} ($dG$-norm of the error). To verify the convergence rates of the PolyDG solution $u_h$ in \eqref{eq:dgPoisson} we use two different scripts: \texttt{RunhConvergenceLaplacian.m} and \texttt{RunpConvergenceLaplacian.m} accounting for the $h$- convergence (mesh size) and the $\ell$-convergence (polynomial degree) respectively.
This time, we set up the data in the external script \\ \noindent \texttt{InputData/DataConvTestLap.m}. Concerning the previous one, four mesh with decreasing granularity $h$ are provided in input within the field \texttt{Data.meshfileseq}.
As the output of the aforementioned scripts we obtain the plot in Figure~\ref{fig:convergence_laplacian}. In particular, on the left, we can observe the convergence of the PolyDG solution obtained with $\ell_\kappa=4$ for any $\kappa \in \mathcal{T}_h$ for the $L^2$- and $dG$-norms, confirming the theoretical results in \eqref{eq:error_estimates}.
On the right, the exponential convergence with respect to the polynomial degree $\ell$ is also shown, cf. \eqref{eq:error_estimates}, by fixing $N_{el} = 100$.
For assessing the performance of \lymph{} in terms of computational time, we consider the Poisson problem solved on a Cartesian grid made of $N = 19600$ elements using a polynomial degree $\ell = 5$. These choices of discretization parameters lead to $411,600$ degrees of freedom. 
Table~\ref{tab:computationalcost_laplacian} presents the computational times, averaged over five simulations. Specifically, we compare the performance of the matrix assembly phase using the quadrature-free and subtriangulation strategies for numerical integration.
We observe that using the quadrature-free approach saves approximately \red{$20\%$} of computational time during the matrix assembly phase. Additionally, this approach becomes increasingly advantageous for general meshes, as it becomes more efficient than the subtriangulation strategy as the number of vertices in the polygonal elements increases, which is  the case with agglomerated grids. Table~\ref{tab:computationalcost_laplacian} also reports the averaged computational times for assembling the right-hand side and solving the corresponding algebraic system. It is important to note that for the right-hand side of the problem, the quadrature-free approach cannot be used as-is, because we need to integrate functions that, in general, are not polynomials. Extending the quadrature-free approach to non-polynomial functions is an open question and will be the subject of further investigation.
The numerical simulations have been performed as a \emph{serial} job, using the \texttt{Kami} cluster (40 computing nodes configured as follows: \textbf{CPU} 2x AMD EPYC 7413 24-Core Processor, 
\textbf{RAM} 512Gb) at the Department of Mathematics, Politecnico di Milano.
\begin{table}[h!]
    \color{black}
    \centering
    \footnotesize
    \begin{tabular}{c|c|c|c|c|c}
        Mat. assembly (QF) &  Mat. assembly (ST) & RHS assembly & Linear system  & \texttt{.csv}-file saving & \texttt{.vtk}-file saving \\ \hline
        25.6500~s  & 31.7024~s & 5.8172~s & 24.7120~s & 9.3620~s & 13.2213~s 
    \end{tabular}
    \caption{\color{black} Test case of Section~\ref{sec:poisson}. Computational times for a test case with $411,600$ \texttt{dofs} ($N = 19600$, $\ell = 5$). We compare the quadrature-free (QF) and subtriangulation (ST) strategies for the numerical evaluation of integrals during matrix assembly. We also report the computational time for the assembly of the right-hand side (RHS) as well as for solving the algebraic system. For the output, we consider two different formats: \texttt{.csv} and \texttt{.vtk} files.}
    \label{tab:computationalcost_laplacian}
\end{table}

\color{black}
\subsubsection{Advantages of polygonal meshes}\label{sec:Poisson_poly_vs_tria}
\begin{figure}
    \centering
    \includegraphics[width=0.24\textwidth]{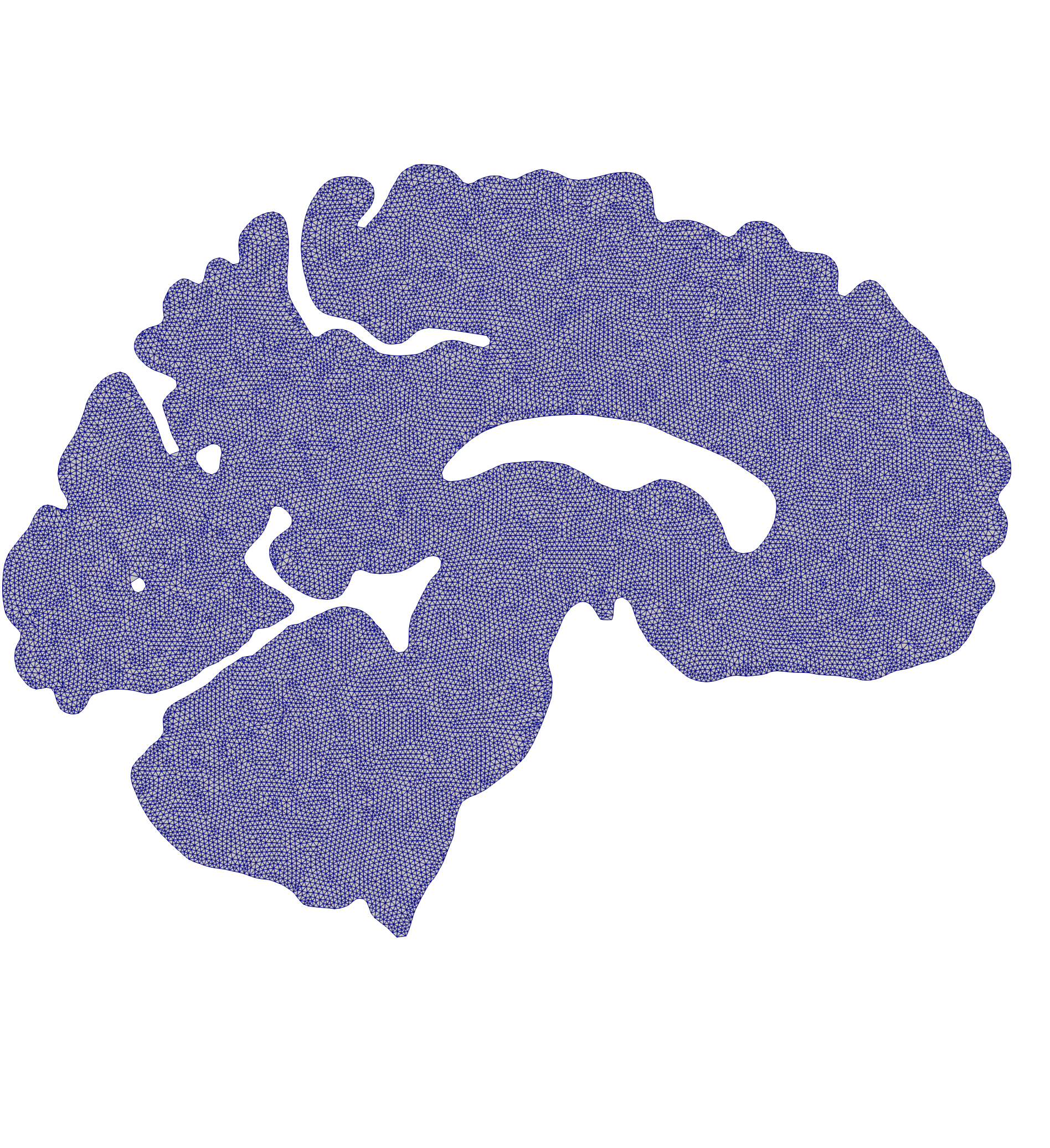}
    \includegraphics[width=0.24\textwidth]{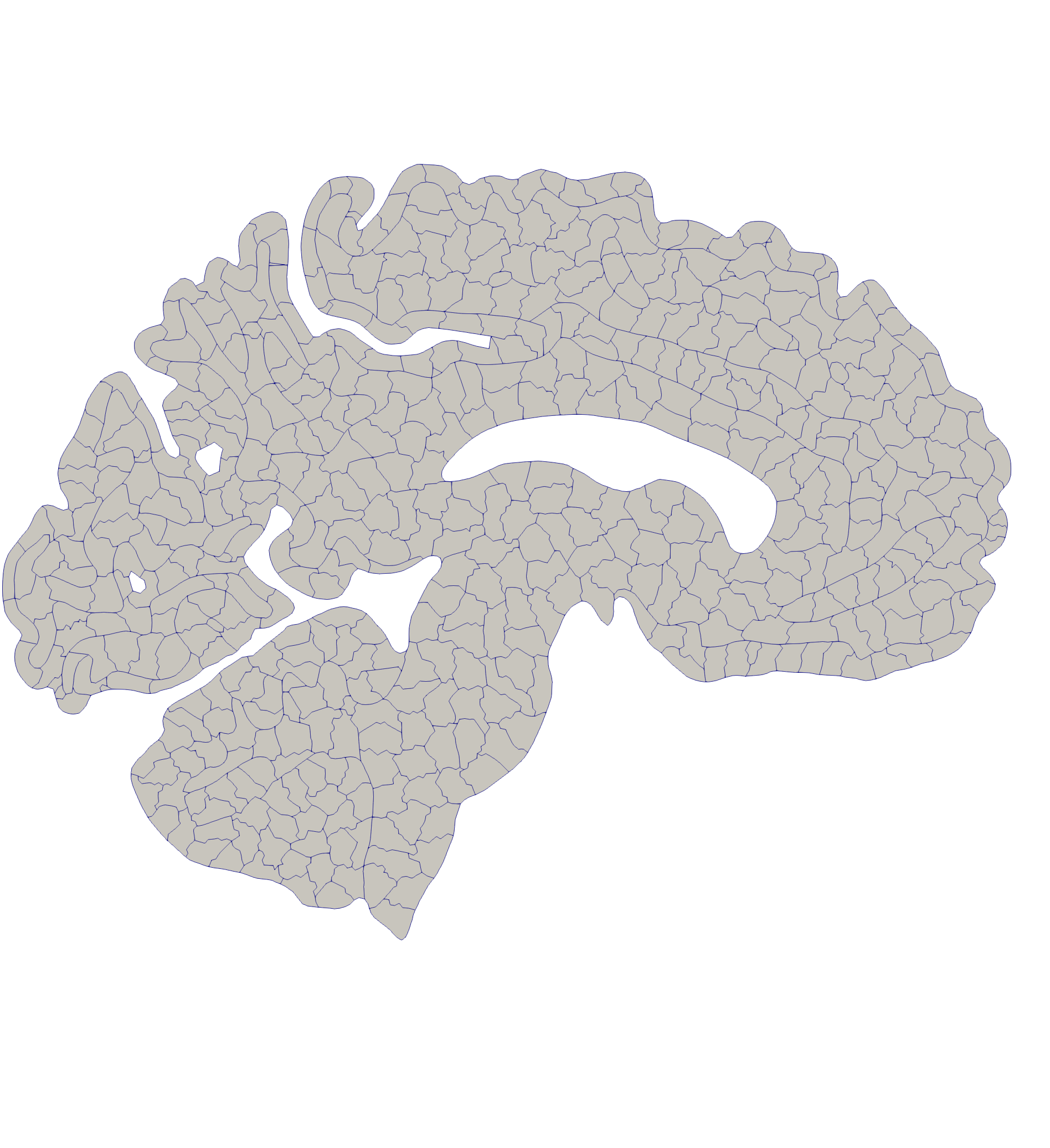}
    \includegraphics[width=0.49\textwidth]{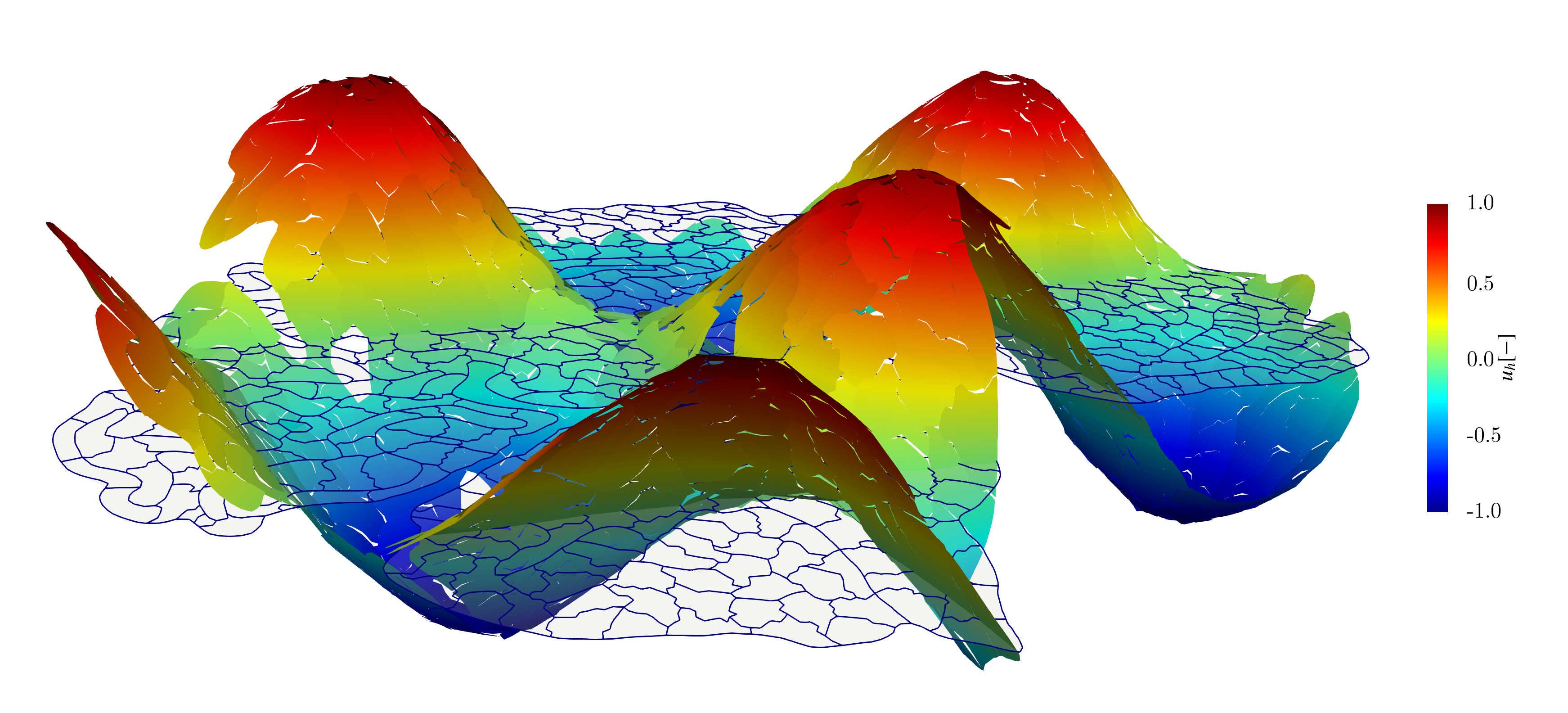}\\
    \caption{Test case of Section~\ref{sec:poisson}. Triangular mesh of brain section (left), agglomerated polygonal mesh (center), and computed PolyDG solution on the polygonal grid with $\ell=1$ (right).}
    \label{fig:sol_lap_comp}
\end{figure}

As an application of the method, we compare the solution of the Poisson problem on a complicated domain, to show the advantages of the proposed solver on general (agglomerated) polygons. In particular, we consider problem \eqref{eq:poisson} with the following data: $\mu(\boldsymbol{x}) = 1$,  $f(\boldsymbol{x}) = 400\pi^2\sin(20\pi x)\cos(20\pi y)$, and $g(\boldsymbol{x}) = \sin(20\pi x)\cos(20\pi y)$, whose exact solution is $u(\boldsymbol{x}) = \sin(20\pi x)\cos(20\pi y)$. The domain $\Omega$ has been derived from a magnetic resonance image of the brain (from the OASIS-3 database \cite{OASIS3}) and it is firstly meshed with $42\,891$ triangular elements (see Figure~\ref{fig:sol_lap_comp}-left). Then we agglomerate the mesh into $534$ polygons (see Figure~\ref{fig:sol_lap_comp}-center), which can capture all the geometric features of the geometry, but with a reduced number of elements in the underlying computational mesh.\\
\par
\begin{figure}
    \centering
    \includegraphics[width=\textwidth]{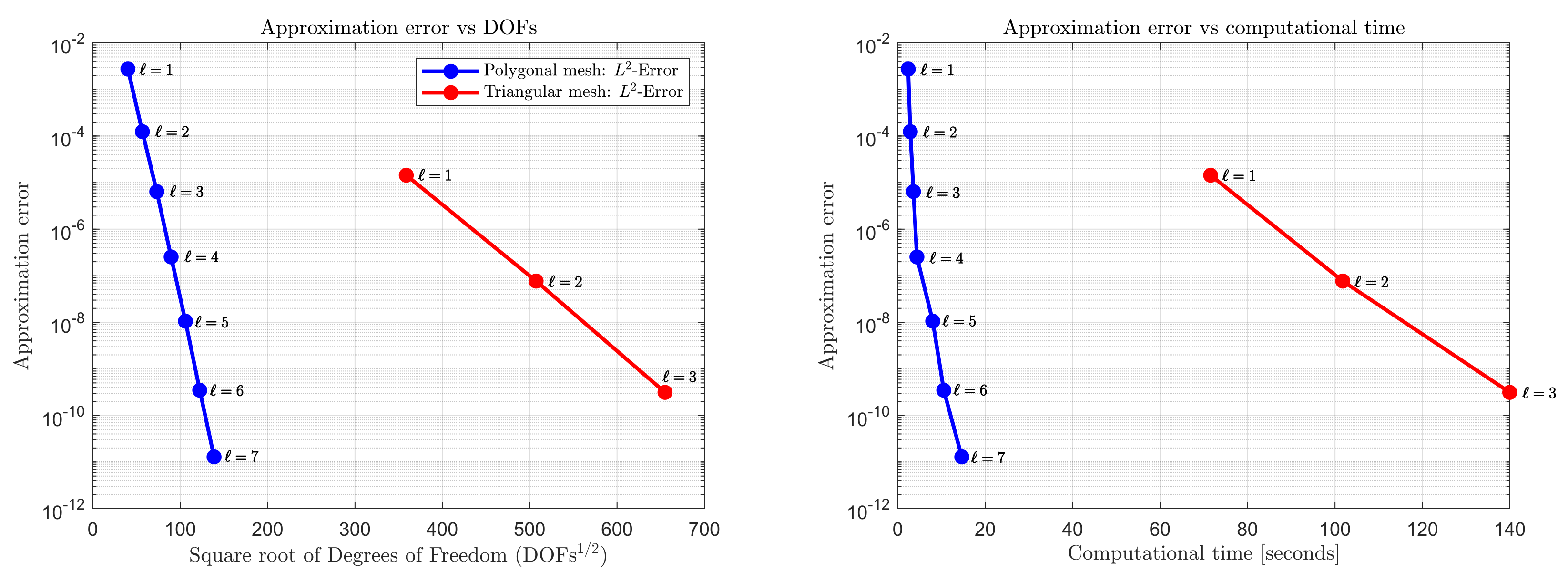}\\
    \caption{\red{Test case of Section~\ref{sec:poisson}.  
    Computed approximation errors in the $L^2$-norm on the triangular (red) and agglomerated (blue) meshes in Figure \ref{fig:sol_lap_comp}. Computed errors versus the total number of degrees of freedom (left), and versus the computational time (right).}}
    \label{fig:err_lap_comp}
\end{figure}
\par
We solve the problem using the triangular mesh with different polynomial degrees $\ell=1,2,3$ and the agglomerated one with $\ell=1,...,7$. In Figure~\ref{fig:sol_lap_comp} (right), we show the solution obtained on the agglomerated grid with $\ell=1$. Moreover, we compare the performance of the solution algorithm, by visualizing the approximation errors versus the number of DOFs (see Figure~\ref{fig:err_lap_comp} left) \red{and versus the computational time including the  matrix and forcing term assembly and linear system solution (see Figure~\ref{fig:err_lap_comp} right)}. 
By fixing the level of accuracy, 
the approximation on the agglomerated grid outperforms the one on the triangular grid. Generally, using a high-order polynomial on an agglomerated mesh is more advantageous than using a low-order approximation on the underlying triangular grid, at least whenever the exact solution is sufficiently regular.


\subsection{The heat equation}
\label{sec:time_dependent_problems}
Let us consider a polygonal domain $\Omega \subset \mathbb{R}^2$, we denote by $\Gamma = \partial \Omega$ its boundary with an outward normal
unit vector $\bm n$. On the boundary, we assume to impose Dirichlet boundary conditions ($ u  = g $). \red{Given a sufficiently regular external load $f$ and initial data $u_0$, the heat equation in $\Omega \times (0,T]$ is given by}
\color{black}
\begin{equation}
    \label{eq:heat}
    \frac{\partial u}{\partial t} - \nabla \cdot (\mu \nabla u) = f, \quad   {\rm in}  \; \Omega \times (0,T],
\end{equation}
\color{black}
with initial condition $u=u_0$, in  $\Omega \times \{ 0 \}$. The PolyDG formulation \cite{Cangiani2017} reads: for any time $t\in (0,T]$ find $ u_h = u_h(t) \in V_h^\ell$ such that 
\color{black}
\begin{equation}\label{eq:heat_dg}
    \sum_{\kappa \in  \mathcal{T}_h} (\dot{u}_h, v_h)_{\kappa} + a_{dG}(u_h, v_h) = \sum_{\kappa \in  \mathcal{T}_h} ( f, v_h)_{\kappa}  \quad \forall v_h \in V_h^\ell,
\end{equation}
\color{black}
with initial conditions $u_h = u_{0h}$, where $u_{0,h}$ is the $L^2$-projection of the initial data on $V_h^\ell$. In \eqref{eq:heat_dg} the bilinear form $a_{dG}(\cdot, \cdot)$ is defined as in the previous section. Now, by introducing a set of basis functions $\{\varphi_j\}_{j=1}^{N_h}$ for $V_h^\ell$ we can easily get the following system of first-order differential equations:
\color{black}
\begin{equation}\label{eq:ode_heat}
  \begin{cases}
      \bm M \dot{U}_h(t)  + A \bm U_h(t) = \bm F(t) & 
      t\in(0,T], \\
      \bm U_h(0) = \bm U_{0h}, & 
  \end{cases}  
\end{equation}
\color{black}
To integrate system \eqref{eq:ode_heat} in time we apply the $\theta-$method scheme. The numerical simulation is performed using the Crank-Nicolson scheme ($\theta=1/2$), see, e.g., \cite{CrankNicolson_1996,Butcher2016}.
We remark that it is possible to couple the \lymph\, library  with time integration schemes already present in MATLAB (e.g., Runge-Kutta schemes \cite{Butcher2016}).

\subsubsection*{Test case with discontinuous boundary conditions}\label{sec:heat_pb}
\begin{figure}
    \centering
    \includegraphics[width=0.49\textwidth]{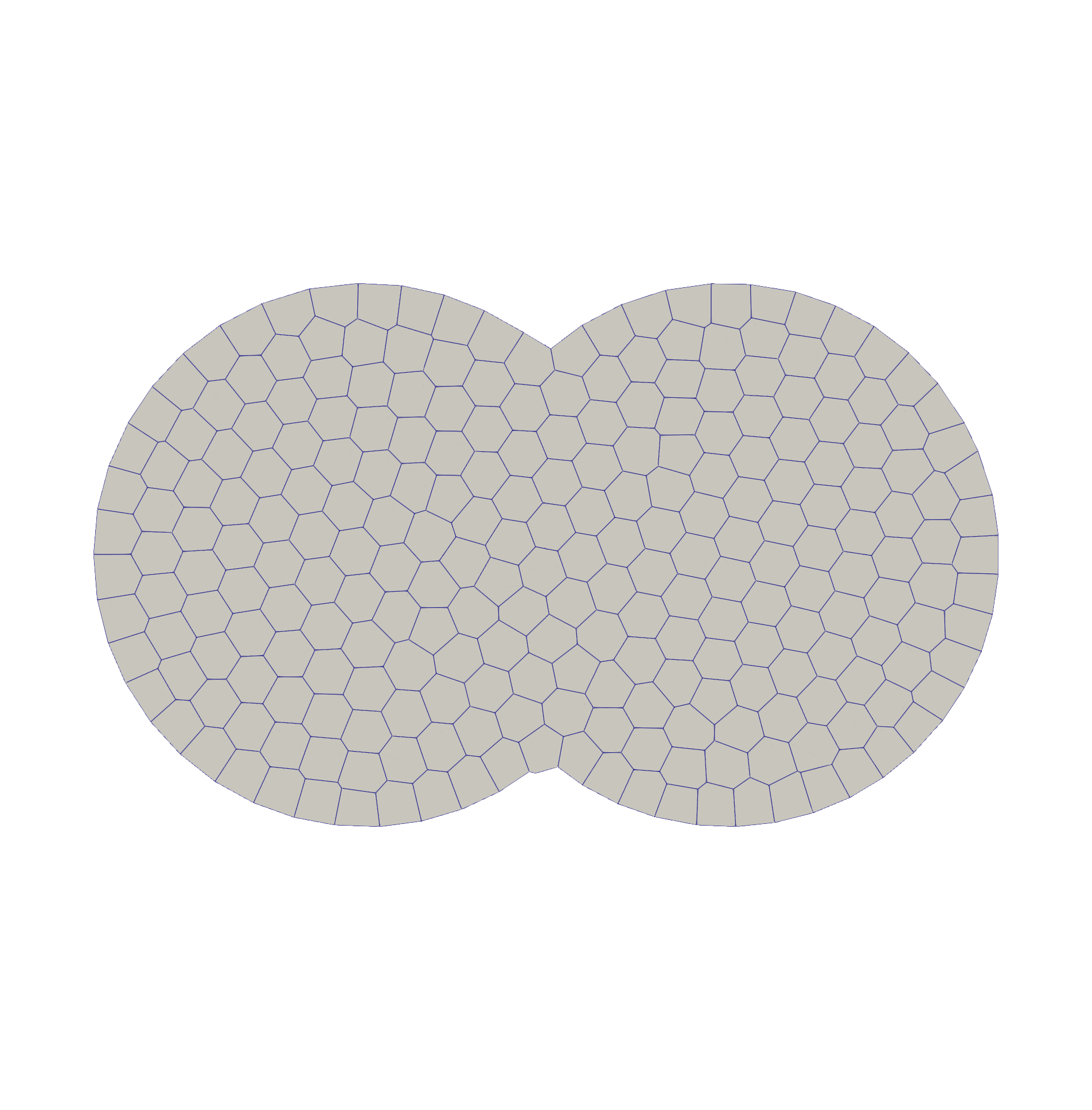}
    \includegraphics[width=0.49\textwidth]{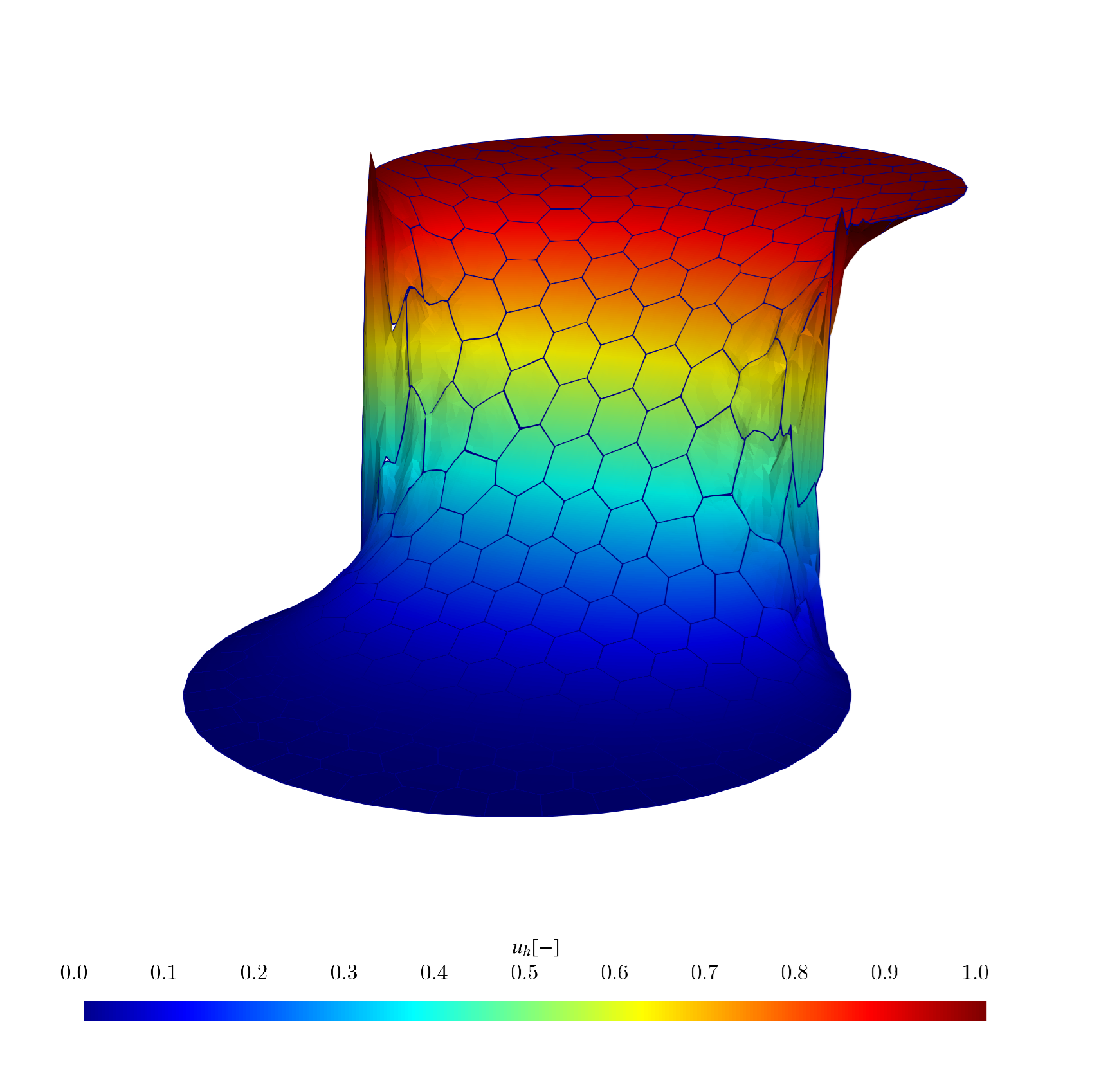}\\
    \caption{Test case of Section~\ref{sec:time_dependent_problems}. Polygonal mesh of the two circles (left) and computed PolyDG solution at time $t=1.0$.}
    \label{fig:sol_heat}
\end{figure}

As an application of the presented PolyDG method, we solve with \lymph~ the heat equation problem presented in \cite{quarteroniEDP}. \red{It considers the parabolic problem with $\mu=0.1$ and a homogeneous forcing term $f=0$}. The domain $\Omega$ is composed of two overlapping circles of radius 0.5 and center $(-0.5,0)$ and $(0.5,0)$ respectively. For the numerical discretization of this problem, we construct a polygonal mesh through the PolyMesher software \cite{talischi2012polymesher} (see Figure~\ref{fig:sol_heat} left). The mesh we adopt in the simulation is composed of 250 elements. 
\par
We consider discontinuous Dirichlet boundary conditions on the top of the model such that $u(\boldsymbol{x})=0$ for $x\leq0$ and $u(\boldsymbol{x})=1$ for $x>0$, where $\boldsymbol{x}=(x,y)$. As shown in the test case reference \cite{quarteroniEDP}, we obtain that the solution at time $t=1$ is smoothed inside the computational domain, due to the diffusion process. The discretization in space is performed employing polynomials of degree $\ell=5$, while the time discretization uses a timestep $\Delta t = 0.002$, and a final time $T=1$.
In Figure~\ref{fig:sol_heat} we report the snapshot of the solution at the final time $T=1$. The solution is coherent to what is expected from the literature \cite{quarteroniEDP}. 

\FloatBarrier
\subsection{The elastodynamics system}\label{sec::elastodynamics}
We consider a bounded convex polygonal domain $\Omega \subset \mathbb{R}^2$, we denote by $\Gamma = \partial \Omega$ its boundary with outward normal
unit vector $\bm n$. The boundary is assumed to be composed of two disjoint portions $\Gamma_D$ and $\Gamma_N$, where Dirichlet ($\bm u  = \bm 0$), and Neumann ($\bm \sigma(\bm u) \bm n = \bm 0$),  conditions are imposed, respectively. 
Given sufficiently regular external loads $\bm f$ and $\bm g$ and initial data $\bm u_0$ and $\bm v_0$, cf. 
 \cite{AntoniettiMazzieri2018},  the equations of (linear) elastodynamics in $\Omega \times (0,T]$ are given by
\begin{equation}\label{eq:elastodynamics}
    \rho \frac{\partial^2\bm u}{\partial t^2} - \nabla \cdot \bm \sigma (\bm u) = \bm f, \quad   {\rm in}  \; \Omega \times (0,T],
\end{equation}
with initial conditions $(\bm u, \frac{\partial\bm u}{\partial t})(0) = (\bm u_0, \bm v_0)$,  in  $\Omega$. We  denote by $\bm u:  \Omega \times [0, T ] \rightarrow \mathbb{R}^2$ the displacement vector and by $\bm \sigma : \Omega \times [0, T ] \rightarrow \mathbb{S}$ the stress tensor where 
$\mathbb {S}$ is the space of symmetric, $2 \times 2$, real-valued tensor fields. We assume the \red{an isotropic-linear elasticity constitutive model} $\bm \sigma (\bm u) =  2\mu \bm \epsilon(\bm u) + \lambda tr(\bm \epsilon(\bm u) )I $, where $\bm \epsilon(\bm u)$ is the symmetric gradient of $\bm u$, $I$ is the identity tensor, $tr(\cdot)$ is the
trace operator, and $\lambda, \mu \in  L^\infty(\Omega)$ are Lamé's parameters. The compressional (P) and shear (S) wave velocities of
the medium are obtained through the relations $c_P =\sqrt{(\lambda + 2\mu)/\rho}$  and $c_S =\sqrt{\mu/\rho}$, respectively.
By following \cite{AntoniettiMazzieri2018}, we obtain the PolyDG formulation: for any time $t\in (0,T]$ find $\bm u_h = \bm u_h(t) \in \bm V_h^\ell$ such that 
\begin{equation}\label{eq:elastodynamics_dg}
    \sum_{\kappa \in  \mathcal{T}_h} (\rho \bm u, \bm v)_{\kappa} + a^e_{dG}(\bm u, \bm v) = \sum_{\kappa \in  \mathcal{T}_h} ( \bm f, \bm v)_{\kappa}  \quad \forall \bm v \in \bm V_h^\ell,
\end{equation}
with initial conditions $(\bm u_h, \frac{\partial\bm u_h}{\partial t}) = (\bm u_{0h}, \bm v_{0h})$, where $\bm u_{0,h}$ and $\bm v_{0,h}$ are the $L^2$-projection of the initial data on $\bm V_h^\ell$. 
In \eqref{eq:elastodynamics_dg} the bilinear forms $a_{dG}^e(\cdot, \cdot)$ is defined  as
\begin{equation}
    a_{dG}^e(\bm u, \bm v)  = \sum_{\kappa \in  \mathcal{T}_h} ( \bm \sigma(\bm u), \bm \epsilon(\bm v))_{\kappa} - \hspace{-3mm} \sum_{e \in \mathcal{F}_h \setminus \Gamma_N} \hspace{-3mm} \Big( \big(\avg{\bm \sigma(\bm u)},\jump{\bm v} \big)_e + \big(\avg{\bm \sigma(\bm v)},\jump{\bm u} \big)_e - \big(\eta_e\jump{(\bm u)},\jump{\bm v} \big)_e \Big) \quad \forall  \bm u, \bm v \in \bm V_h^\ell,
\end{equation}
with $\eta_e$ as in \cite[eq. (9)]{AntoniettiMazzieri2018}. Now, by introducing a set of basis functions $\{\bm \varphi^1_j, \bm \varphi^2_j\}_{j=1}^{N_h}$ for $\bm V_h^\ell$ we can easily get the following system of second order differential equations:
\begin{equation}\label{eq:ode_elastodynamics}
  \begin{cases}
      M \ddot{\bm U}_h(t)  + A \bm U_h(t) = \bm F(t) & 
      t\in(0,T], \\
      \red{\dot{\bm U}_h(0)=\bm V_{0h}, \qquad \bm U_h(0) = \bm U_{0h}}, & 
  \end{cases}  
\end{equation}
To integrate system \eqref{eq:ode_elastodynamics}
in time we apply the Newmark $\beta$-scheme, with $\beta=\frac14$ and $\gamma=\frac12$. 

\begin{figure}
    \centering
    \includegraphics[width=0.7\textwidth]{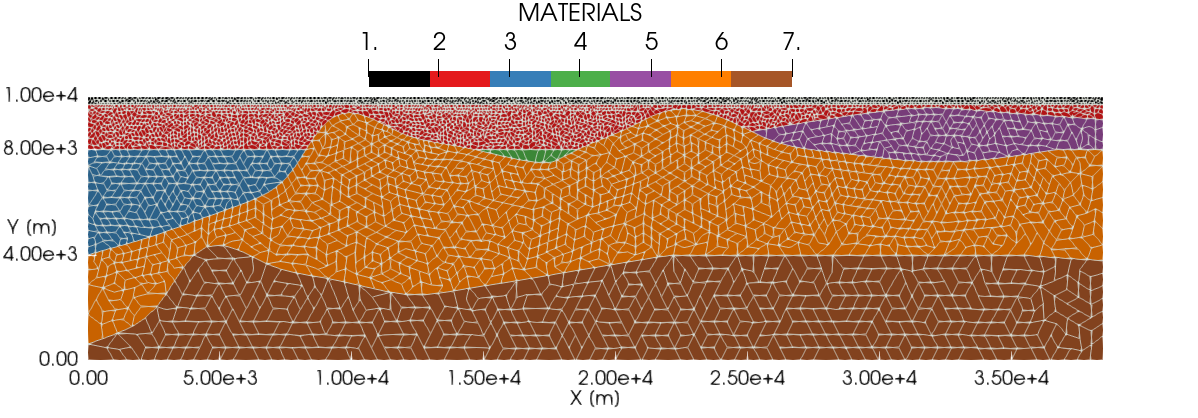}
    \caption{Test case of Section~\ref{sec::elastodynamics}. Unstructured polygonal grid with mesh spacing of about $h \approx 160$ m for material 1 to $h \approx 1500$ m for material 7; cf. Table~\ref{tab:tab_materials}. Mesh file is available in \texttt{Elastodynamics/InputMeshPhysics/MeshEmilia.mat.}}\label{fig:mesh_elasto}
\end{figure}

As an application of the presented PolyDG method, we solve with \lymph~ the wave propagation problem presented in  \cite[Section 5.4.3.2]{Antonietti2021Review}. It considers the elastic wave propagation $\Omega = (0, 38.4)$ km $\times (0, 10)$ km representing an idealized bidimensional Earth’s cross-section, see Figure~\ref{fig:mesh_elasto}. 
\begin{table}
    \centering
    \footnotesize
    \begin{tabular}{c|c|c|c|c|c|c|c}
       Materials  & 1 & 2 & 3 & 4 & 5 & 6 & 7 \\
       \hline
        $\rho$ [kg/m$^3$]  & 1800 & 1800 & 2050 & 2050 & 2050 & 2400 & 2450 \\
        \hline
        $c_S$ [m/s] & 294 & 450 & 600 & 600 & 600 & 1515 & 1600 \\ \hline 
        $c_P$ [m/s] & 1321 & 2024 & 1920 & 1920 & 1920 & 3030 & 3200\\
    \end{tabular}
    \caption{Test case of Section~\ref{sec::elastodynamics}. Material properties used for the computational domain in Figure~\ref{fig:mesh_elasto}, cf. also \texttt{Elastodynamics/InputData/Elastic/DataTestPhysicsEla.m}}\label{tab:tab_materials}
\end{table}

\begin{figure}
    \centering
    \includegraphics[width=0.49\textwidth]{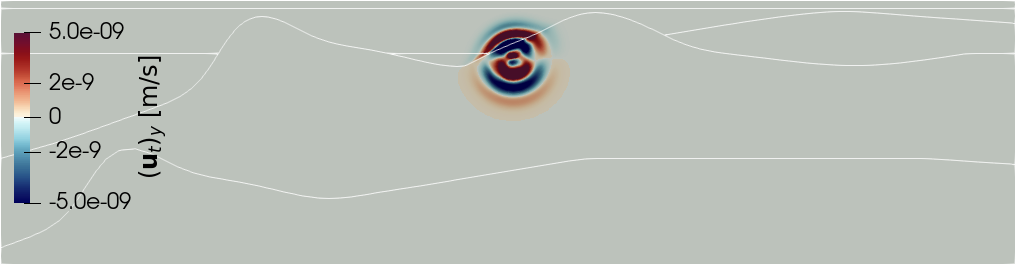}
    \includegraphics[width=0.49\textwidth]{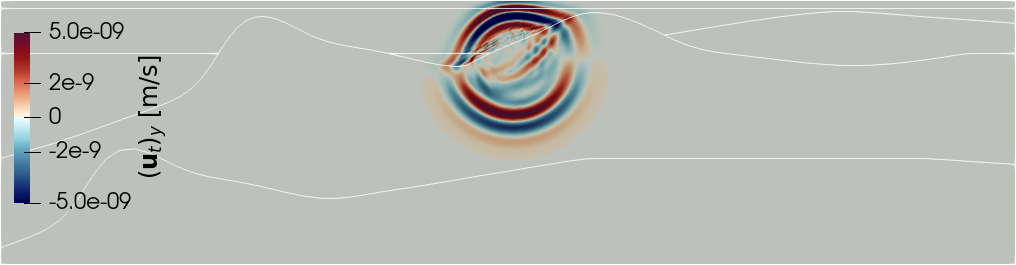}\\
    \includegraphics[width=0.49\textwidth]{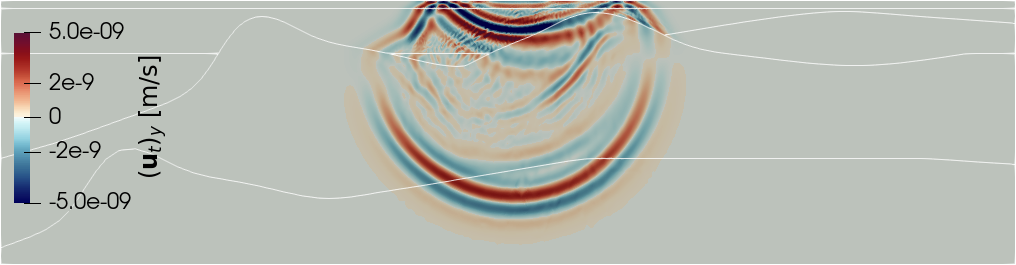}
    \includegraphics[width=0.49\textwidth]{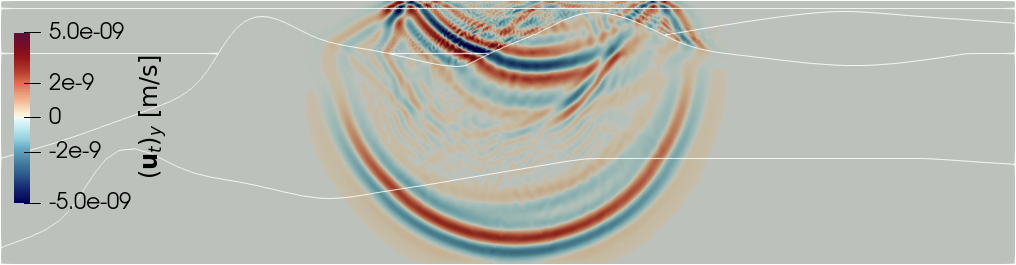}
    \caption{Test case of Section~\ref{sec::elastodynamics}.  Snapshots of the computed vertical velocity $(\bm u_t)_y$ at different times $t = 0.75$ (top-left), $t=1.25$ (top-right), $t=2.25$ (bottom-left), $t=2.75$ (bottom-right). Due to the material heterogeneity's, high oscillations and perturbations of the wavefront can be observed, as well as the effect of the free surface on top.}
    \label{fig:snapshots_wave}
\end{figure}

We consider homogeneous Neumann conditions on the top of the model ($\bm \sigma \bm n = \bm 0)$ whereas homogeneous Dirichlet conditions ($\bm u = \bm 0$) are set on the remaining boundaries. The bottom and the lateral boundaries are set far enough from the point source to prevent any reﬂections from the boundaries of the waves of interest. 
We simulate a double-couple moment source load of the form
$\bm f(\bm x, t) = -  I \cdot \nabla \delta (\bm x - \bm x_s)S(t)$,
where $\delta(\bm x - \bm x_s)$ is the Dirac delta distribution centered in $\bm x_s = (19432,7800)$~m and $S(t) = (1-8\pi^2(t-0.5)^2)e^{-4\pi^2(t-0.5)^2}$ is the source time function. We assign constant material properties within each region as described in Table \ref{tab:tab_materials}. The computational domain is discretized using an unstructured grid consisting of 4870 (agglomerated) polygonal elements, with a mesh size varying from $h \approx 160$ m for
material $1$ to $h \approx 1500$ m for material $7$; cf. Table \ref{tab:tab_materials}. 
We consider also a polynomial degree $\ell=5$, $\Delta t = 0.001$~s, and a final time $T=4$~s.
In Figure~\ref{fig:snapshots_wave} we report a set of snapshots of the computed vertical velocity field $(\bm u_t)_y$. The discontinuities between the mechanical properties of the materials produce oscillations and perturbations on the wavefront; surface waves are visible. 
Finally, we report in Table~\ref{tab:comput_time} the computing time of the presented test (averaged over $5$ different simulations). The numerical simulation has been performed as a \emph{serial} job, using the \texttt{Kami} cluster described before.

\begin{table}[h!]
    \centering
    \footnotesize
    \color{black}
    \begin{tabular}{c|c|c|c|c|c}
        Mat. assembly (QF) &  Mat. assembly (ST) & RHS assembly & Linear system  & \texttt{.csv}-file saving & \texttt{.vtk}-file saving \\ \hline
        10.9928~s & 13.0421~s & 2.3805~s & 62.4727~s & 3.4967~s & 4.3340~s 
    \end{tabular}
    \caption{\color{black} Test case of Section~\ref{sec::elastodynamics}. 
    Computational times for a test case with $409,080$ \texttt{dofs} ($N = 4870$, $\ell = 5$). We compare the quadrature-free (QF) and subtriangulation (ST) strategies for the numerical evaluation of integrals during matrix assembly. We also report the computational time for the assembly of the right-hand side (RHS) as well as for solving the algebraic system. For the output, we consider two different formats: \texttt{.csv} and \texttt{.vtk} files.}
    \label{tab:comput_time}
\end{table}
\color{black}

\FloatBarrier
\subsection{A multi-physics problem: the poroelasto-acoustic system}\label{sec:coupled-poro-acoustic}
In this last application, we present a wave propagation problem in a coupled poroelastic-acoustic domain $\Omega = \Omega_p \cup \Omega_a$, where
$\Omega_p$ and $\Omega_a$ are the poroelastic and acoustic medium, respectively. These domains share an interface denoted as $\Gamma_I$. \\

In the fluid domain $\Omega_a$ we consider the following model
\begin{equation}\label{eq::acoustic}
\begin{cases}
c^{-2}\dfrac{\partial^2 \varphi_a}{\partial t^2} - \rho_a^{-1}  \nabla\cdot(\rho_a \nabla \varphi_a)=f_a,  & \text{in }\Omega_a\times(0,T],\\
  \varphi_a = g_a, & \text{on }\Gamma_{aD}\times(0,T],\\ 
  \nabla \varphi_a \cdot \bm n_a = 0, & \text{on }\Gamma_{aN}\times(0,T], \\
  \red{\varphi_a(0) = \varphi_0, \qquad \dfrac{\partial \varphi_a}{\partial t}(0)=\psi_0,} & \text{in }\Omega_a,
\end{cases}
\end{equation}
where $\varphi_a$ is the acoustic potential, $c>0$ is the wave speed, $\rho_a >0$ is the medium density, and $f_a, g_a, \varphi_0$ and $\psi_0$ are regular enough data. We consider the boundary $\Gamma_a = \Gamma_{aD}\cup \Gamma_{aN}$ of $\Omega_a$ to be sufficiently regular, with the normal $\bm n_a$ outward pointing, and to be decomposed into two distinct portions denoted by $\Gamma_{aD}$ (resp. $\Gamma_{aN})$ where Dirichlet (resp. Neumann) conditions are imposed. \\

In the porous domain $\Omega_p$ we consider the system 
\cite{biot1}:
\begin{equation}\label{eq::poroelasticity}
\begin{cases}
\rho_p \dfrac{\partial^2 \bm u_p}{\partial t^2} + 
\rho_f\dfrac{\partial^2 \bm u_f}{\partial t^2}
-\nabla\cdot\bm{\sigma}_p(\bm{u}_p,\bm{u}_f)=\bm{f}_p,  &\text{in }\Omega_p\times(0,T],\\[2mm]
\rho_f\dfrac{\partial^2 \bm u_p}{\partial t^2} + 
\rho_w\dfrac{\partial^2 \bm u_f}{\partial t^2} + 
\frac{\eta}{k}\dfrac{\partial \bm u_f}{\partial t}
+\nabla p_p(\bm{u}_p,\bm{u}_f)=\bm{f}_f,  &\text{in }\Omega_p\times(0,T],\\
\bm{\sigma}_p(\bm{u}_p,\bm{u}_f) = \bm{\sigma}_e(\bm{u}_p) 
-\beta\, p_p(\bm{u}_p,\bm{u}_f) I, &\text{in }\Omega_p\times(0,T],\\
 p_p(\bm{u}_p,\bm{u}_f) = -m(\beta \nabla\cdot\bm{u}_p+\nabla\cdot \bm{u}_f),
&\text{in }\Omega_p\times(0,T],\\
(\bm u_p,\bm u_f)  = (\bm g_p, \bm g_f), & \text{on }\Gamma_p\times(0,T],\\[2mm]
  \red{\bm u_p(0) = \bm u_{p0}, \qquad \dfrac{\partial \bm u_p}{\partial t}(0) = \bm v_{p0},} & \text{in }\Omega_p, \\[2mm]
  \red{\bm u_f(0) = \bm u_{f0}, \qquad \dfrac{\partial \bm u_f}{\partial t}(0) = \bm v_{f0},} & \text{in }\Omega_p.
\end{cases}
\end{equation}
Here, $\bm{u}_p$ and $\bm{u}_f$ represent the displacements of the solids and the filtration, respectively, $\rho_p$ is the average density given by $\rho_p=\phi\rho_f+(1-\phi)\rho_s$, where $\rho_s>0$ is the solid density, $\rho_f>0$ is the density of the saturated fluid, $\rho_w$ is defined as $\rho_w=\frac{a}{\phi}\rho_f$, with porosity $\phi$ satisfying $0<\phi_0\leq\phi\leq\phi_1<1$, and the tortuosity $a>1$ measures the deviation of the fluid paths from straight streamlines. The dynamic viscosity of the fluid is given by $\eta>0$, the absolute permeability by $k>0$, while the Biot--Willis coefficient $\beta$ and the Biot modulus $m$ are such that $\phi<\beta\le1$ and $m\ge m_0>0$. In  \eqref{eq::poroelasticity}, $\bm{f}_p,\bm{f}_f,\bm g_p,\bm g_f, \bm u_{p0}, \bm v_{p0}, \bm u_{f0}$, and $\bm v_{f0}$ are given (regular enough) data.
Moreover, the elastic tensor $\bm{\sigma}_e(\bm{u}_p)$ is defined by Hooke's law introduced in Section \ref{sec::elastodynamics}. Finally, we suppose the boundary $\Gamma_p$ of $\Omega_p$ to be sufficiently regular, having outward pointing unit normal $\bm n_p$. \\

On the shared interface $\Gamma_I$, the following coupling conditions are prescribed  
\begin{equation}\begin{cases} \label{eq:interface}
-\bm{\sigma}_p(\bm u_p,\bm u_f) \bm{n}_p  = \rho_a\dfrac{\partial \varphi_a}{\partial t}\bm{n}_p, & \textrm{ on }\Gamma_{I}  \times (0,T],    \\
p_p(\bm u_p , \bm u_f ) = \rho_a \dfrac{\partial \varphi_a}{\partial t},
 & \textrm{ on }\Gamma_{I}  \times (0,T],  \\
\red{-\left(\dfrac{\partial \bm{u}_p}{\partial t}+\dfrac{\partial \bm{u}_f}{\partial t}\right)\cdot\bm{n}_p = \nabla\varphi_a\cdot\bm{n}_p,} & \textrm{ on }\Gamma_{I}  \times (0,T],
\end{cases}\end{equation}
expressing the continuity of normal stresses, continuity of pressure, and conservation of mass, respectively. We refer the reader to \cite{Antonietti2021} for a detailed description of the problem as well its numerical discretization and analysis. 
For this problem the \textit{polytopic} mesh $\mathcal{T}_h=\mathcal{T}^p_h\cup\mathcal{T}^a_h$, where $\mathcal{T}^{\#}_h=\{\kappa\in\mathcal{T}_h:\overline{\kappa}\subseteq\overline{\Omega}_{\#}\}$, with $\#=\{p,a\}$.
Implicit in this decomposition is the assumption that the meshes $\mathcal{T}_h^a$ and $\mathcal{T}_h^p$ are aligned with $\Omega_a$ and $\Omega_p$, respectively. 
Moreover, we denote by $V_{\#h}^\ell = {V_h^\ell}_{|_{\Omega_{\#}}}$,  with $\#=\{p,a\}$.
As explained in \cite{Antonietti2021}, the space discretization of \eqref{eq::acoustic}-\eqref{eq::poroelasticity}-\eqref{eq:interface} with a PolydG method leads to the following problem:
for $t\in(0,T]$, find $(\bm{u}_{ph},\bm{u}_{fh},\varphi_{ah})(t)\in \bm{V}_{ph}^\ell\times \bm{V}_{ph}^\ell\times V_{ah}^\ell$, s.t.
\begin{multline}
 \sum_{\kappa \in  \mathcal{T}_{ph}} 
 (\rho \ddot{\bm{u}}_{ph} + \rho_f \ddot{\bm{u}}_{fh} , \bm{v}_h)_{\kappa} 
    +  ( \rho_f \ddot{\bm{u}}_{ph} + \rho_w \ddot{\bm{u}}_{fh}, \bm z_h)_{\kappa} + 
    \sum_{\kappa \in  \mathcal{T}_{ah}} (c^{-2} \ddot{\varphi}_{ah}, \psi_h)_{\kappa} +   
  \sum_{\kappa \in  \mathcal{T}_{ph}} (\eta k^{-1} \dot{\bm{u}}_{fh},\bm{z}_h)_{\kappa} 
    \\
+ \mathcal{A}_h^e(\bm{u}_{ph},\bm{v}_h)+
\mathcal{A}_h^p(\beta\bm{u}_{ph}+\bm {u}_{fh},\beta\bm{v}_h+\bm{z}_h) +  \mathcal{A}_h^a(\varphi_{ah},\psi_h),
\\ + \sum_{e \in \Gamma_I} (\rho_a \dot{\varphi}_{ah} , (\bm{v}_h+\bm{z}_h)\cdot\bm{n}_p )_e
- \sum_{e \in \Gamma_I} ((\dot{\bm{u}}_{ph}+\dot{\bm{u}}_{fh})\cdot\bm{n}_p, \rho_a \psi_{h} )_e
=  \sum_{\kappa \in  \mathcal{T}_{ph}} (\bm{f}_p,\bm{v}_h)_\kappa + (\bm{f}_f,\bm{z}_h)_\kappa + \sum_{\kappa \in  \mathcal{T}_{ah}} (f_a,\psi_h)_\kappa \label{eq::dgsystem}
\end{multline} 
for any $(\bm{v}_h,\bm{z}_h,\psi_h) \in \bm{V}_{ph}^\ell\times \bm{V}_{ph}^\ell\times V_{ah}^\ell$. Note that in \eqref{eq::dgsystem}
 we have used the notation $\dot{u}$ (resp. $\ddot{u}$) to indicate the first (resp. second) time derivative of a generic function $u$. 
The bilinear forms $\mathcal{A}_h^e(\cdot,\cdot),
\mathcal{A}_h^p(\cdot,\cdot)$, and $  \mathcal{A}_h^a(\cdot,\cdot)$ appearing in \eqref{eq::dgsystem} are defined as in \cite[eq. (3.4)]{Antonietti2021}.
Now, by introducing a set of basis functions for $\bm V_{ph}^\ell$ and $V_{ah}^\ell$ and denoting by ($U_{pf}$, $U_{fh}$,$\Phi_{ah}$) the vector of the expansion coefficients in the chosen basis of the unknowns $\bm{u}_{ph}$,  $\bm{u}_{fh}$ and $\varphi_{ah}$, respectively, we obtain the following system of ordinary differential equations:
\begin{multline}\label{eq::algebraic}
\left[ \begin{matrix} 
\bm{M}_\rho^p & \bm{M}_{\rho_f}^p & 0 \\
\bm{M}_{\rho_f}^p & \bm{M}_{\rho_w}^p & 0 \\
0 & 0 & \bm{M}_{c^{-2}}^a
\end{matrix} \right]  
\left[ \begin{matrix} 
\ddot{U}_{ph} \\
\ddot{U}_{fh} \\ 
\ddot{\Phi}_{ah}
\end{matrix} \right] + 
\left[ \begin{matrix} 
0 & 0 & \bm{C}^p \\
0 & \bm{B} & \bm{C}^p \\
\bm{C}^a & \bm{C}^a & 0
\end{matrix} \right]  
\left[ \begin{matrix} 
\dot{U}_{ph} \\
\dot{U}_{fh} \\ 
\dot{\Phi}_{ah}
\end{matrix} \right] 
+\left[\begin{matrix} 
\bm{A}^e +  \bm{A}_{\beta^2}^p &  \bm{A}_{\beta}^p & 0 \\
\bm{A}_{\beta}^p & \bm{A}^p & 0 \\ 
0 & 0 & \bm{A}^a
\end{matrix} \right] 
\left[ \begin{matrix} 
U_{ph} \\
U_{fh} \\ 
\Phi_{ah}
\end{matrix} \right] = 
\left[ \begin{matrix} 
\bm{F}^p \\
\bm{F}^f \\ 
\bm{F}^a
\end{matrix} \right] 
\end{multline}
with initial conditions $(U_{ph},\dot{U}_{ph})(0) =(U_{p0},V_{p0})$,  $(U_{fh},\dot{U}_{fh})(0) = (U_{f0},V_{f0})$, and $(\Phi_{ah},\dot{\Phi}_{ah}) (0)=(\Phi_{0},\Psi_0)$. In the above system  $\bm{F}^p$, $\bm{F}^f$ and $\bm{F}^a$ are the vector representations of the right-hand side of \eqref{eq::dgsystem}.
To integrate in time system \eqref{eq::algebraic} we use the Newmark scheme as in Section \ref{sec::elastodynamics}. We refer the interested reader to \cite{Antonietti2021} for all the details.\\

\red{We point out that the assembly strategy of \lymph{} allows for system \eqref{eq::algebraic} to be assembled with almost no additional complexity with respect to assembling each of the involved physics separately. Indeed, the\\ \noindent \texttt{Physics/PoroElastoAcoustics/Assembly} folder contains a sub-folder with the routines to assemble the terms corresponding to each of the involved physics, an additional one to assemble the coupling terms, and a function to combine all the terms in the monolithic matrices.}

We now consider the problem of acoustic scattering by a porous cylinder.
\begin{figure}
    \centering    \includegraphics[width=0.55\textwidth]{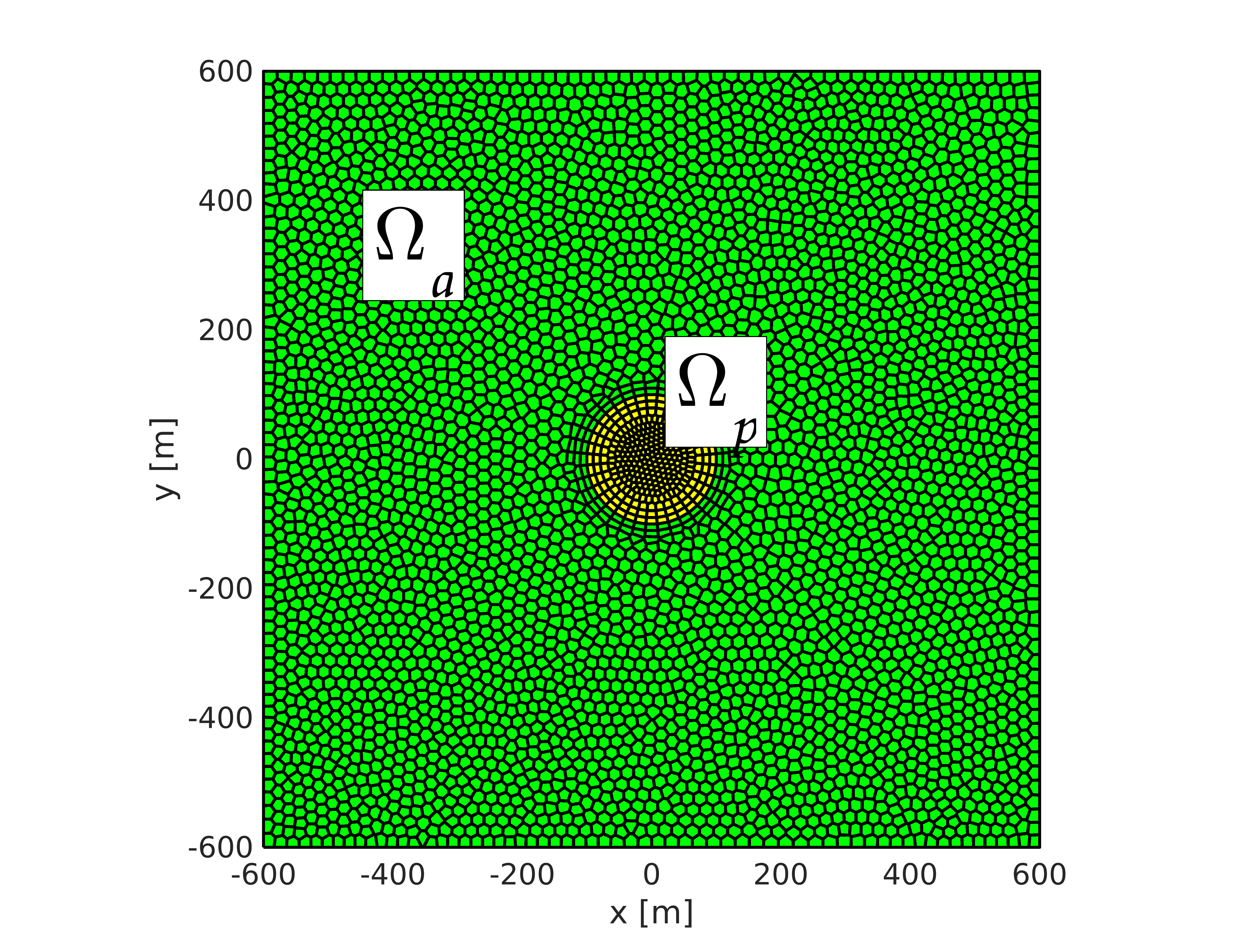}    \includegraphics[width=0.44\textwidth]{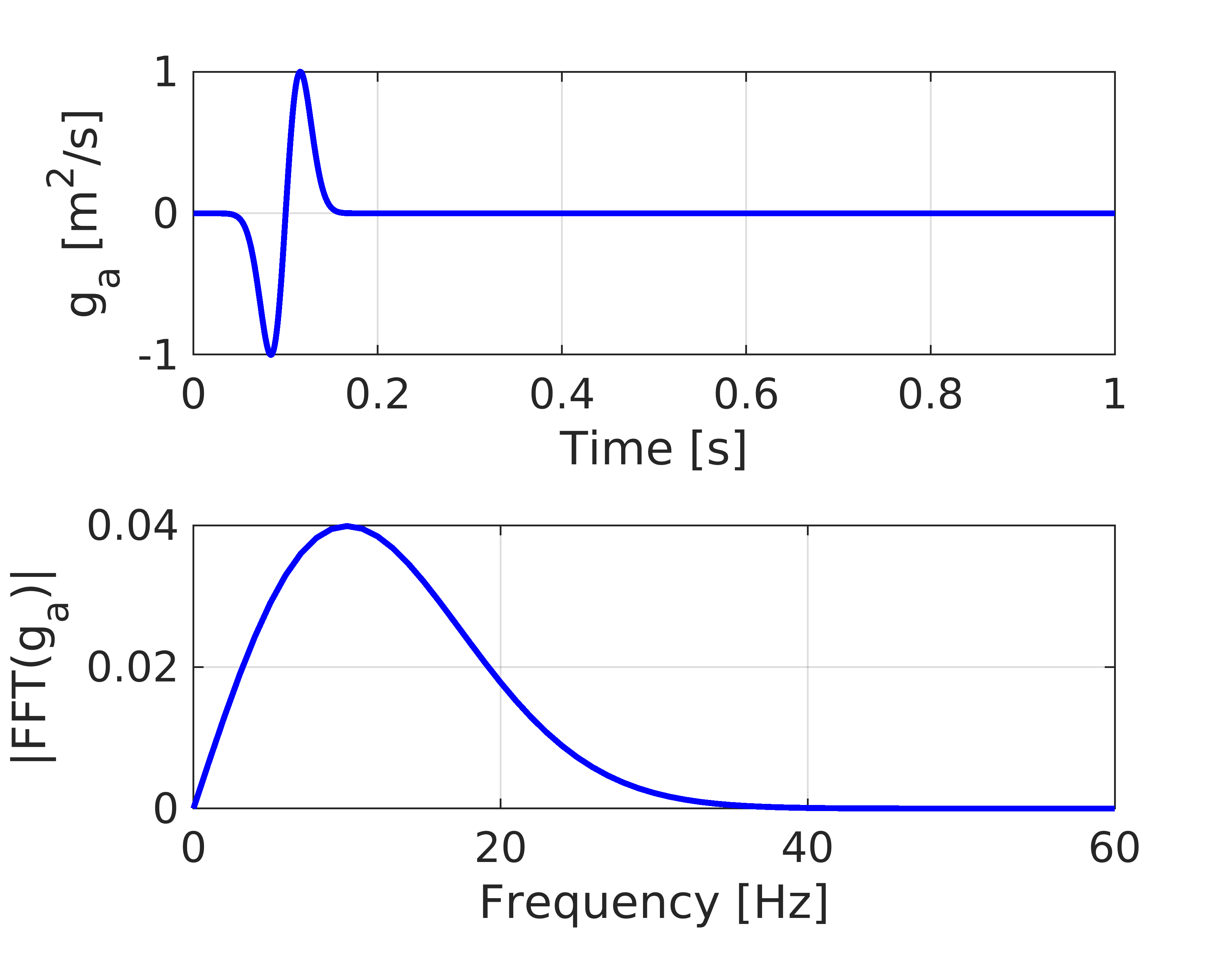}
    \caption{Test case of Section~\ref{sec:coupled-poro-acoustic}. Left: circular porous domain $\Omega_p$ (yellow) surrounded by an acoustic medium $\Omega_a$ (green). The mesh is made by $4287$ polygonal elements divided into $3979$ for $\Omega_a$ ($h\approx 36.5~m$) and $308$ for $\Omega_p$ ($h\approx 11.8~m$). Right: Dirichlet boundary condition applied on the domain's bottom edge (top), the absolute value of its Fourier transform (bottom).}
    \label{fig:AcousticCavity-domain}
\end{figure}
The domain is given by $\Omega = (-600,600)~m^2$, see Figure~\ref{fig:AcousticCavity-domain}, cf. also \cite{chiavassa_lombard_2013}, where
a circular porous cylinder $\Omega_p$ is surrounded by the acoustic medium $\Omega_a$. More precisely, we consider a circular interface $\Gamma_I$ of radius $100~m$ centered at $(0, 0)$, dividing $\Omega$ into a circular porous domain $\Omega_p$ (yellow, in Figure~\ref{fig:AcousticCavity-domain}) surrounded by an acoustic medium $\Omega_a$ (green, in Figure~\ref{fig:AcousticCavity-domain}).    
The computational domain is discretized by 4287 polygonal elements having a characteristic size $h\approx 36~m$, and we consider a space polynomial degree equal to 4. A locally refined mesh is employed around the interface to compute precisely the various wave conversions, cf. Figure~\ref{fig:AcousticCavity-domain}. The source is an acoustic plane wave given as a Dirichlet condition on the bottom boundary, i.e., $\Gamma_{aD} = 
 (-600,600)\times \{-600\}$, cf. Figure ~\ref{fig:AcousticCavity-domain} (right).
The pulse is a Gaussian wavelet 
$g_a(t) = 2\pi f_p \sqrt{e}(t - 1/f_p)e^{-2(\pi f_p)^2 (t - 1/f_p)^2}$,
having a peak frequency $f_p = 10~Hz$. 
A sound soft condition is enforced on the remaining boundaries, i.e., $\nabla \varphi_a \cdot \bm n_a = 0$ on $\Gamma_{aN}$. The physical parameters for this test case are listed in Table \ref{tab::table_poroacoustic}.
%
%
\begin{table}[htbp]
\centering
\begin{tabular}{llllll}
\cline{1-6}
\textbf{Fluid}  & Fluid density      & $\rho_f$    & 1000    & $\rm kg/m^3$    &  \\
                &                    & $\rho_a$    & 1000   & $\rm kg/m^3$    &  \\
                & Wave velocity      & $c$         & 1500   & $\rm m/s$       &  \\
                & Dynamic viscosity  & $\eta$      & $1.05\cdot 10^{-3}$      & $\rm Pa\cdot s$ &  \\ \cline{1-6}
\textbf{Grain}  & Solid density      & $\rho_s$    & 2690   & $\rm kg/m^3$    &  \\
                & Shear modulus      & $\mu$       & 1.86 $\cdot 10^9$      & $\rm Pa$   &  \\ \cline{1-6}
\textbf{Matrix} & Porosity           & $\phi$      & 0.38    &             &  \\
                & Tortuosity         & $a$         & 1.8      &             &  \\
& Permeability  & $k$                & $2.79\cdot 10^{-11}$    & $\rm m^2$       &  \\
& Lam\'e coefficient   & $\lambda$ & 1.20$\cdot 10^8$ & $\rm Pa$        &  \\
& Biot's coefficient & $m$         & 5.34$\cdot 10^9$ & $\rm Pa$        &  \\
& Biot's coefficient & $\beta$     & 0.95             &             &  \\  \cline{1-6}
\end{tabular}
\caption{Test case of Section~\ref{sec:coupled-poro-acoustic}.  Physical parameters for the poroelastic-acoustic test case.}
\label{tab::table_poroacoustic}
\end{table}
%
\begin{figure}
    \centering   \includegraphics[width=0.32\textwidth]{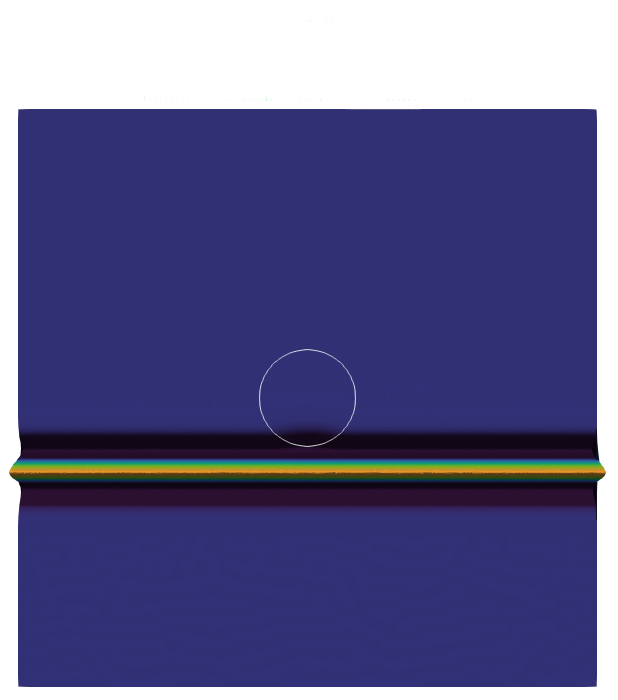} \includegraphics[width=0.32\textwidth]{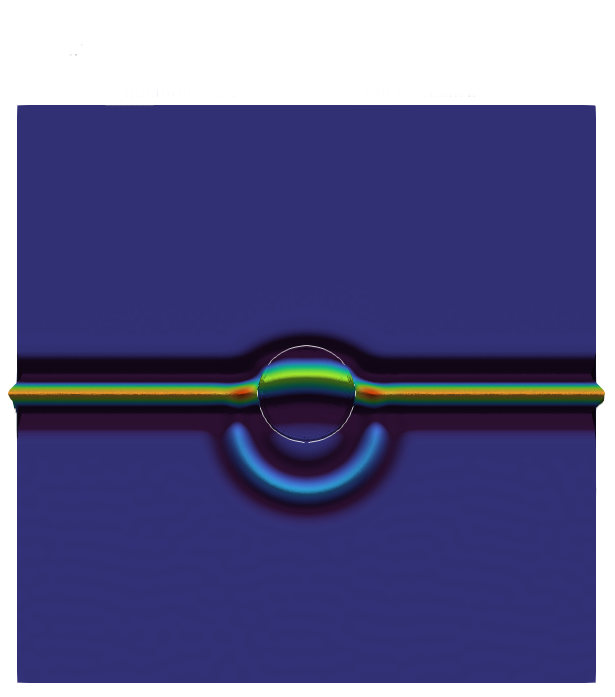} \includegraphics[width=0.32\textwidth]{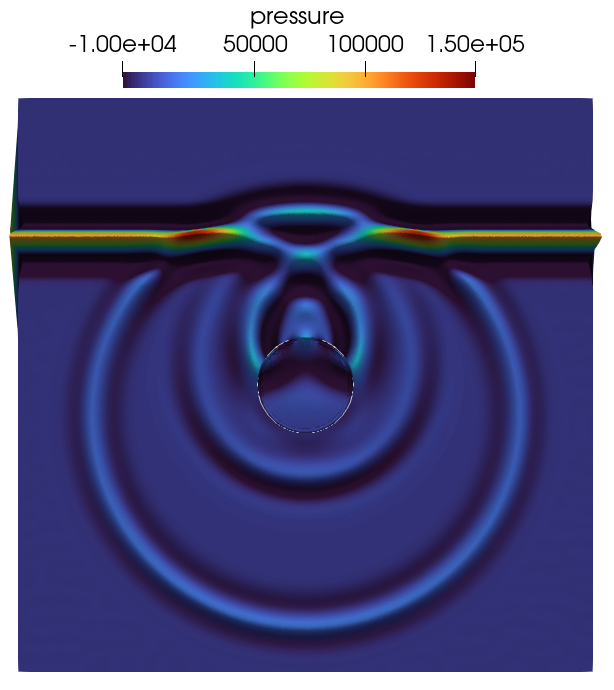} 
    \caption{Test case of Section~\ref{sec:coupled-poro-acoustic}. Snapshots of the computed pressure field at different time instants: $t=0.4~s$ (top-left), $t=0.5~s$ (top-right), $t=0.6~s$ (bottom-left), and $t=0.7~s$ (bottom-right).}    \label{fig:snapshot_acousticcavity}
\end{figure}
In Figure~\ref{fig:snapshot_acousticcavity} we plot the computed pressure field $p_a = \rho_a \dot{\varphi}_a(t)$ in $\Omega_a$ and $p_p = -m(\beta \nabla\cdot\bm{u}_p + \nabla\cdot \bm{u}_f)$ in $\Omega_p$ at different time instants, considering a final time $T=1~s$, and a time step $\Delta t = 0.001~s$.
It is possible to observe that the wave is moving from the bottom to the top of the domain, impacting the porous cylinder, and producing a scattered circular wavefield directed backwards. The numerical results are aligned with the corresponding ones presented in \cite{chiavassa_lombard_2013}.

\section{Conclusions}
This paper presents \lymph{}, a general-purpose Matlab library for the approximate solution of multi-physics differential problems in two-dimensions. For the spatial discretization of the underlying differential systems,  \lymph{} library is based on high-order discontinuous Galerkin methods on polytopal grids, making its use attractive for several areas of engineering and applied sciences applications. The target of this paper is to introduce the library step-by-step and to show the potential of the software, starting from the solution of classical differential problems.  As \lymph{} is a user-friendly, general-purpose library, the authors think that its use can be widely extended to other engineering applications. Interesting future developments of this work include the design of more robust and flexible (agglomeration-driven) mesh generation algorithms and the introduction of $hp$-refinement approaches, and the extension to three-dimension.

\section{Acknowledgements}
\red{We thank the Editors Lin-bo Zhang and Wolfgang Bangerth, and the anonymous Reviewers for their valuable insights and comments, which greatly  contributed to improving the content of the present work.}
We acknowledge Prof. Paul Houston and Dr. Giorgio Pennesi for the original implementation of the quadrature-free approach \cite{AntoniettiHoustonPennesi_18} adopted in the library.
The brain MRI images were provided by OASIS-3: Longitudinal Multimodal Neuroimaging: Principal Investigators: T. Benzinger, D. Marcus, J. Morris; NIH P30 AG066444, P50 AG00561, P30 NS09857781, P01 AG026276, P01 AG003991, R01 AG043434, UL1 TR000448, R01 EB009352. AV-45 doses were provided by Avid Radiopharmaceuticals, a wholly-owned subsidiary of Eli Lilly.
This work received funding from the European Union (ERC SyG, NEMESIS, project number 101115663). Views and opinions expressed are however those of the authors only and do not necessarily reflect those of the European Union or the European Research Council Executive Agency. Neither the European Union nor the granting authority can be held responsible for them. 
PFA, IF, IM, have been partially supported by ICSC--Centro Nazionale di Ricerca in High Performance Computing, Big Data, and Quantum Computing funded by European Union--NextGenerationEU.
PFA and SB have been partially funded by MUR, PRIN 2020 research grant n. 20204LN5N5.
All the authors are members of INdAM-GNCS. The work of IM has been partially supported by the INdAM-GNCS project CUP E53C22001930001.

\bibliographystyle{ieeetr}
\bibliography{bibliography_rev.bib}

\end{document}